\documentclass{conm-p-l-version}

\usepackage{xr}
\newcommand{\bfc}{{\Bbb C}}
\newcommand{\bfz}{{\Bbb Z}} 
\newcommand{\calc}{{\mathcal C}}
\newcommand{\co}{{\mathcal O}}
\newcommand{\cald}{{\mathcal D}}
\newcommand{\calv}{{\mathcal V}}
\newcommand{\cala}{{\mathcal A}}
\newcommand{\calw}{{\mathcal W}}
\newcommand{\cale}{{\mathcal E}}
\newcommand{\calb}{{\mathcal B}}
\newcommand{\cals}{{\mathcal S}}
\newcommand{\csj}{\C_{\Sigma_j}}
\newcommand{\csi}{\C_{\Sigma_i}}
\newcommand{\gj}{\G_{\C,j}^2}
\newcommand{\bnu}{\mathbf{\nu}}
\newcommand{\xst}{X_k^{st}}  
\newcommand{\nxst}{\widetilde{X_k^{st}}}  
\newcommand{\G}{\Gamma}
\newcommand{\C}{{\calc}} 
\newcommand{\noi}{\noindent}
\newcommand{\gxk}{\G(X_{k})}
\newcommand{\gxg}{\G(X_{k},g)}
\newcommand{\gc}{\G_\C}
\newcommand{\gce}{\G_{\C}^1}
\newcommand{\gck}{\G_{\C}^2}
\newcommand{\gcce}{G_{\C}^1}

\newcommand{\fcd}{F_{c,d}}
\newcommand{\Wedge}{{\mathsf W_{\eta, M}}}
\newcommand{\de}{\mathbf{D}}    
\newcommand{\dec}{\mathbf{D}_c} 
\newcommand{\ded}{\mathbf{D}_d} 
\newcommand{\deo}{\mathbf{D}_0} 
\newcommand{\dee}{\mathbf{D}^{(1)}} 
\newcommand{\dek}{\mathbf{D}^{(2)}}  
\newcommand{\deh}{\mathbf{D}^{(3)}}  
\newcommand{\wP}{\widetilde{\Phi}}
\newcommand{\wfcd}{\widetilde{F_{c,d}}}
\newcommand{\ic}{irreducible component}
\newcommand{\alahuz}[1]{$\underline{\mbox{#1}}$}
\newcommand{\bs}{\bigskip}
\newcommand{\sm}{\smallskip}

\def\fig#1{
 \begin{tabular}{c}
    \includegraphics{#1}
  \end{tabular}
}

\begin{document}

\title[Resolution Graphs of Iomdin Series]{Resolution Graphs of Some Surface Singularities, II.\\
(Generalized Iomdin Series)}
\author{Andr\'as N\'emethi}
\address{Department of Mathematics, The Ohio State University, Columbus,
Ohio 43210}
\email{nemethi@math.ohio-state.edu}
\author{\'Agnes  Szil\'ard}
\address{Department of Mathematics, The Ohio State University, Columbus,
Ohio 43210}
\email{szilard@math.ohio-state.edu}
\subjclass{Primary 32S50; Secondary 32S05, 32S25, 14J17, 14B05}
\keywords{surface singularitiy, topological series of singularities, 
resolution graph,  non-isolated singularity,
isolated complete intersection singularity}

\begin{abstract}
In this article, we determine the resolution graph of the
hypersurface singularity $(\{f+g^k=0\},0)$, where $(f,g):({\bfc}^3,0)\to
({\bfc}^2,0)$ is an ICIS, and $k$ is a sufficiently large integer.
All these graphs  are coordinated by an ``universal bi-colored graph''
$\G_{\C}$ associated with the ICIS $(f,g)$. Its definition is rather
involved, and in concrete examples it is difficult to compute it.
Nevertheless, we  present a large number of examples. This is 
very helpful in the exemplification  of its properties as well. 
Then we present our main construction which provides
the dual resolution graph of the series $\{f+g^k=0\}$ from the  graph
$\G_{\C}$  and the integer $k$. This is formulated in a purely combinatorial
algorithm. The result is a highly non-trivial generalization of the
cyclic covering case of $({\bfc}^2,0)$.
\end{abstract}

\maketitle

\noi{\sc Introduction.}

\vspace{2mm}

The present article (in the sequel: Part II) is a natural continuation of 
\cite{partI} (called Part I) to such an extent that we use and cite the 
notations and results of Part I without reviewing them. (Note that
  the sections are numbered continuously through Part I and II.) 

Consider an isolated complete intersection
singularity (ICIS) $\Phi=(f,g):(\bfc^3,0)\to (\bfc^2,0)$ such that $f$ has an 
1--dimensional singular locus. Then construct the 
so called ``generalized Iomdin series'' $f+g^k$ (for $k\geq k_0$)
(see Section \ref{sec4}, i.e.  the first section of Part II). 
Our goal is to describe  the resolution graph $\gxk$ 
of  $(X_k,0):=(\{f+g^k=0\},0)$. 

The Iomdin series is a generalization of the cyclic covering case $f(x,y)+z^k$,
where $f$ is an isolated plane curve singularity. But the differences
(both technical and of principle) are huge. For example, while in the case 
of the cyclic coverings, the relevant monodromy representations are
easy (being  representations of ${\bfz}$), the monodromy  representation
of an ICIS is incomparably more complicated. Moreover, in the cyclic case 
there is  a global Galois action, which disappears in the general case.

However, the construction of $\gxk$ will show strong similarities with the
case of cyclic coverings: we will construct from the ICIS $\Phi$   a universal
graph $\gc$ which will coordinate  the whole 
series of graphs $\{\gxk\}_k$. Already the
existence of such a graph is really remarkable. 

The construction of $\gc$ is rather involved. We will show that all the 
geometric information needed for the construction of the graphs $\{\gxk\}_k$
is contained in a small neighborhood of a special curve arrangement $\C$,
situated in the embedded resolution of the local divisor $(\{fg=0\},0)\subset
(\bfc^3,0)$. $\gc$ is exactly the dual  graph of $\C$. We will decorate this
graph with different weights, which codify the topology of the irreducible
components of $\C$, and also the local behavior of the functions
$f$ and $g$ near $\C$. In Section \ref{sec5} 
we present this construction, some 
examples and also  some properties of $\gc$. E.g., we will prove that
from $\gc$ one can recover the resolution graph of $(\{f=0\},0)$, and also 
the equisingular type of the transversal singularities of $Sing\{f=0\}$. 

In Section \ref{sec6} 
we  present the algorithm which  provides, in a purely 
combinatorial way, the graph $\gxk$ from $\gc$ and the integer $k$ (for
$k\geq k_0$). The validity of the algorithm is based on the following
important fact. The graph $\gxk$ appears as a cyclic covering graph 
of $\gc$ (modified with some Hirzebruch--Jung strings). The weights of
$\gc$ codify all the local data from a neighborhood of $\C$, but a priori
it is not clear that from this 
data one can recover the necessary global information
in order to identify the covering $\gxk$. 
But, using the classification results of Section \ref{sec1} (Part I), 
and the properties of $\gc$, we show that the covering data of the covering 
can be determined from the weights of $\gc$, and with the determined 
covering data there is a unique  covering graph. This basically guarantees
the vanishing of the global invariants. 

The graph $\gc$ can also be connected  with some other invariants of the 
Iomdin series (for details, see (\ref{invs}) and (\ref{fr})). 

\section{Preliminaries about generalized Iomdin series}\label{sec4}

\subsection*{Basics of ICIS.}

\vspace{2mm}

\bekezdes{}{icisbasic} \  In this subsection
 we review  the basic notions related to
 isolated complete intersection singularities (ICIS). For details, see
\cite{Lj}. In order to keep the
 notation uniform,  we only present the case $(\bfc^3,0)\to (\bfc^2,0)$. 

 Consider an analytic germ $\Phi=(f,g):(\bfc^3,0)\to (\bfc^2,0)$
which defines an ICIS. This means that if 
 $I\subset\co_{\bfc^3,0}$ denotes  the ideal generated by $f,g$ and
the $2\times 2$ minors of the Jacobian matrix $(d\Phi)$, then 
 $\dim \co_{\bfc^3,0}/I<\infty$. In other words,
the scheme-theoretical intersection $\Phi^{-1}(0)=
\{ f=0\}\cap\{ g=0\}$ has only an isolated singularity at the origin.
In particular,   $\ (Sing\{f=0\})\cap \{g=0\}= \{ 0\}$ and
 $ \{f=0\}$ intersects $\{g=0\}$ in the complement of the
origin transversally in a smooth (punctured) curve.

 The {\em critical locus} $(C_\Phi,0)$  of $\Phi$ is the set of points of 
$(\bfc^3,0)$, where $\Phi$ is not a submersion. Its image 
$(\Delta_{\Phi},0):=\Phi(C_\Phi,0)\subset(\bfc^2,0)$ is called the
{\em discriminant locus}  of $\Phi$. 
In the sequel, we denote  by $(c,d)$ the local coordinates of $(\bfc^2,0)$.

\bekezdes{The Milnor fibration.}{icisfib}
 Fix a germ  $\Phi$ as above.
 Then  there exist a sufficiently small closed ball
$B^3_\epsilon\subset\bfc^3$ (centered at the origin with radius 
$\epsilon>0$)  and   $D^2_\eta\subset\bfc^2$ with radius
 $0<\eta<<\epsilon$ such that:

$\bullet$\  the set $(\Phi^{-1}(0)\setminus\{0\})\cap B^3_\epsilon$ 
is non-singular;

$\bullet$\  $\partial B^3_{\epsilon'}$ intersects $\Phi^{-1}(0)$ 
      transversally for all $0<\epsilon'\leq\epsilon$.

$\bullet$\ $C_\Phi\cap\Phi^{-1}(D^2_\eta)\cap\partial B^3_\epsilon=\emptyset$;
and   the restriction of $\Phi$ to 
      $\Phi^{-1}(D^2_\eta)\cap\partial B^3_\epsilon$ is a submersion.

\bekezdes{Definition.}{goodrepr}
{\em  The map $\Phi:\,\Phi^{-1}(D^2_\eta)\cap B^3_\epsilon\rightarrow
D^2_\eta$ with the above properties is called a ``good'' representative of the
ICIS $\Phi$.  For such a good representative, let 
$\Sigma_\Phi$ denote the intersection $C_\Phi\cap\Phi^{-1}(D^2_\eta)$ and
  $\Delta_\Phi=\Phi(\Sigma_\Phi)$. }

\vspace{2mm}

With these notations we have the following additional properties
(cf. 2.8 in \cite{Lj}):

\bekezdes{Theorem.}{fifibration} {\em 

\noi (i)\  
$\Phi:\overline{B^3_\epsilon}\cap\Phi^{-1}(D^2_\eta)\longrightarrow D^2_\eta$ 
is proper. The analytic sets 
 $\Sigma_\Phi$ and $\Delta_\Phi$ are 1--dimensional, 
and the restriction $\Phi|_{\Sigma_\Phi}:\,\Sigma_\Phi\rightarrow\Delta_\Phi$
is proper with finite fibres.

\noi (ii)\  $\Phi:
(\Phi^{-1}(D^2_\eta-\Delta_\Phi)\cap B^3_\epsilon,
\Phi^{-1}(D^2_\eta-\Delta_\Phi)\cap \partial B^3_\epsilon)
\rightarrow D^2_\eta-\Delta_\Phi$ is a 
locally trivial $C^\infty$-fibration of a pair of spaces.}

\bekezdes{Definition.}{milnorfibre}
{\em A fibre $\fcd=\Phi^{-1}(c,d)\cap B^3_\epsilon$,
 for $(c,d)\in D^2_\eta-\Delta_\Phi$,  is
called ``a Milnor fibre'' and the fibration itself 
is referred to as ``the Milnor fibration''.
For any fixed base point $b_0=(c_0,d_0)\subset D^2_\eta-\Delta_\Phi$,
one has the natural {\em geometric monodromy representation}:
$m_{geom}:\, \pi_1(D^2_\eta-\Delta_\Phi, b_0)\ \longrightarrow\ 
\mbox{\it Diff}\,^\infty (F_{b_0})/isotopy.$}

\vspace{1mm}

In general,  it is extremely difficult to determine this representation.
Actually, even the induced {\em algebraic monodromy representation} 
$m_{alg}:\, \pi_1(D^2_\eta-\Delta_\Phi, b_0)\ \to
 Aut H_*(F_{b_0},{\bfz})$ is a  very complicated object. 

\bekezdes{}{d}\
Let $\Delta_\Phi=\Delta_1\cup...\cup\Delta_t$ denote the decomposition
of  $\Delta_\Phi$  into irreducible components.
Let $\Sigma_{i,1},...,\Sigma_{i,s_i}$ be the irreducible decomposition
of $\Phi^{-1}(\Delta_i)\cap \Sigma_\Phi$ for $i=1,...,t$. Clearly,
$\Sigma_\Phi=\cup_{i=1}^t\cup_{j=1}^{s_i}\Sigma_{ij}$ is a decomposition
of $\Sigma_\Phi$ into irreducible components.

Assume that $f$ has a 1-dimensional singular locus. 
Clearly, $Sing\{f=0\}\subset\Sigma_\Phi$ and
$\Phi(Sing\{f=0\})=\{c=0\}$ is an irreducible component of the discriminant
$\Delta_\Phi$. By convention, we denote this component by $\Delta_1$.
Then $\Sigma_{1,1},...,\Sigma_{1,s_1}$ are exactly the 
irreducible components of $(Sing\{f=0\},0)$. Sometimes, we will use the
simplified notations $\Sigma_{1,j}=\Sigma_j$ and  $s_1=s$. 
Part (i) of the above theorem  guarantees that
the restriction $\Phi: \Sigma_{j}\rightarrow\Delta_1$ is a branched covering
for any $1\leq j\leq s$. Let $d_{j}$ denote its degree.

\bekezdes{Definition. -- Transversal singularity types.}{transversdef} 
Fix an irreducible component $\Sigma_{j}$ of 
$Sing\{f=0\}$, $j\in\{1,...,s\}$.
 Take a point $q\in\Sigma_{j}-\{0\}$ and the germ $(H,q)$ of a generic smooth 
transversal slice to $\Sigma_j$  at $q$. 
The intersection $(\{f=0\}\cap H,q)$ determines  an isolated plane curve
singularity $(T\Sigma_j,q)\subset (H,q)$. Its topological (or equisingular)
type  does not depend on the choice
of $q$ and $(H,q)$. {\em It is called  the 
transversal singularity of  $Sing\{f=0\}$ corresponding
to the branch $\Sigma_{j}$}.

\subsection*{Generalized Iomdin series.}

\vspace{2mm}

\noindent The {\em generalized Iomdin series} is a special family   of 
topological series of composed hypersurface  singularities.
For  the general theory of the topological
series,  the interested reader can consult \cite{Schr,NS,NP,NZ}.
In this subsection we restrict ourselves to the definition of the 
special case of the Iomdin series. 

\bekezdes{Definition.}{goal2} \ 
{\em Consider a germ $f:(\bfc^3,0)\rightarrow(\bfc,0)$ with a 1-dimensional
singular locus. Then fix another germ 
$g:(\bfc^3,0)\rightarrow(\bfc,0)$, such that
the pair $\Phi=(f,g):(\bfc^3,0)\rightarrow(\bfc^2,0)$ forms  an ICIS.
Then the family of singularities  
$f+g^k:(\bfc^3,0)\rightarrow(\bfc,0)$ 
constitutes the {\em generalized Iomdin series} of $f$, associated with $g$,
where $k$ is any integer 
 larger than a bound $k_0=k_{0,\Phi}$. (The bound depends only on $\Phi$,
and will be determined in the sequel). }

\vspace{2mm}

Now, we provide a more precise description in terms of the 
resolution graphs of $f$ and $f+g^k$. First notice that 
the germs  $f$ and $f+g^k$ are composed singularities.
Indeed, $f=\Phi\circ P$, where 
$P:(\bfc^2,0)\rightarrow(\bfc,0)$ is the  plane curve singularity
$(c,d)\mapsto c$. Moreover, $f+g^k=\Phi\circ P_k$ 
with the same $\Phi$, but   $P_k$ given by  $(c,d)\mapsto c+d^k$.
Similarly as above, let $\Delta_{\Phi}$ be the discriminant locus of 
$\Phi$ with the distinguished component $\Delta_1=\{c=0\}$. 
A possible embedded resolution graph of the divisor
$(P^{-1}(0)\cup\Delta_{\Phi},0)\subset (\bfc^2,0)$ 
can be visualized schematically as:

\vspace{2.2cm}\hspace{-1cm}
\fig{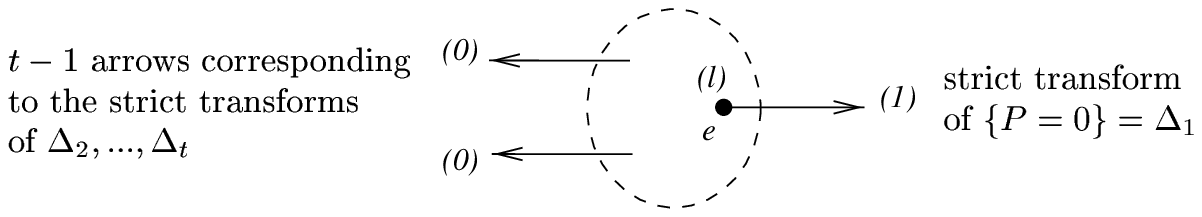}


\noi Above, the multiplicities $(m)$ are the vanishing orders of the (lift 
of the) germ $P$ along the corresponding components of the total transform
of $P^{-1}(0)\cup \Delta_{\Phi}$.

By definition,  the germ $f+g^k$ belongs to the generalized Iomdin series
of $f$, associated with $g$, if $k\geq l+1$. This means that a possible 
(schematic)   embedded resolution graph of $P^{-1}_k(0)\cup\Delta_{\Phi}$ is:

\vspace*{3cm} 

\hspace{-1cm}\fig{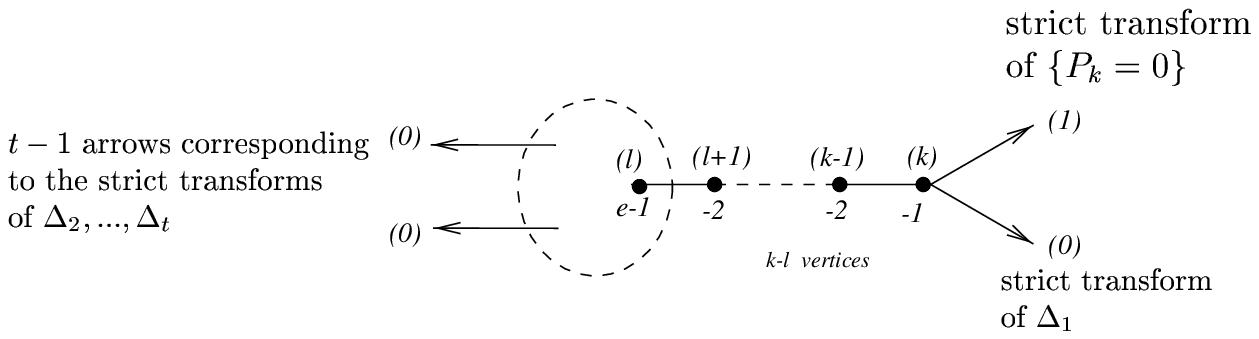}

\noi  where  the multiplicities $(m)$ are the vanishing orders of the germ
$P_k$ along the corresponding components of the total transform
of $P^{-1}_k(0)\cup \Delta_{\Phi}$. (Hence $k_{0,\Phi}=l+1$.)

\bekezdes{Example.}{tn}  Let $f=x^n+y^n+xyz$ with $n\geq 3$ and $g=z$. Then 
$\Phi=(f,g): (\bfc^3,0)\rightarrow (\bfc^2,0)$ is an ICIS.
The critical locus of $\Phi$ is
$\Sigma_\Phi=\{x=y=0\}\cup\,\cup_{j=1}^n
              \{ x=\xi_0^jy,\ z=-n\xi_0^{j(n-1)}y^{n-2}\},$
where $\xi_0=e^{2\pi i/n}$. Moreover, $\{x=y=0\}=Sing(f^{-1}(0))$.
Then $\Phi(\{x=y=0\})=\{c=0\}$ and
$\Phi( \{ x=\xi_0^jy,\ z=-n\xi_0^{j(n-1)}y^{n-2}\})=
\{ c=(2-n)y^n,\ d=-n\xi_0^{j(n-1)}y^{n-2}\}$ (for $j=1,...,n$), 
but some of these components may coincide.  Indeed,
let $l=\gcd(n,n-2)$. Note that $l=1$ or $l=2$. 

If $l=2$, i.e. if  $n=2m$ for some $m\geq 1$, then 
$\Delta_\Phi=\Delta_1\cup\Delta_2\cup\Delta_3$, where
$\Delta_1=\{c=0\},\  \Delta_2=\{c^{m-1}=Kd^m\},\ 
\Delta_3=\{c^{m-1}=-Kd^m\}$ (for some 
constant  $K\neq 0$). 

\noindent
The resolution graph of $P^{-1}(0)\cup\Delta_{\Phi}=
\Delta_1\cup\Delta_2\cup\Delta_3$ is

\vspace*{3cm}

\hspace{-1cm}
\fig{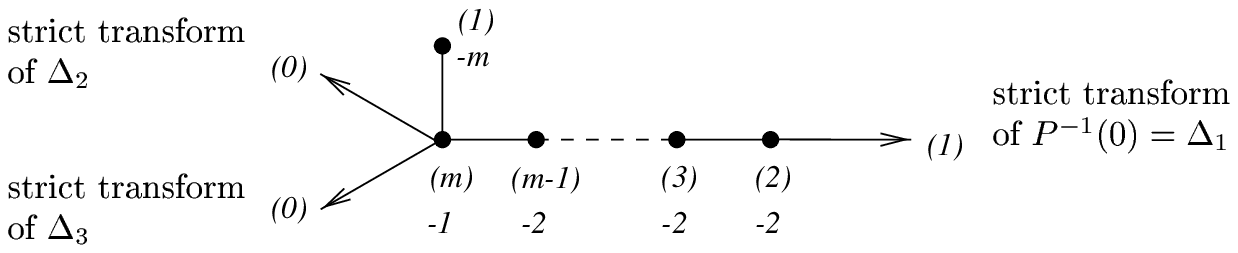}

\noi
 The resolution graph of $P_k^{-1}(0)\cup\Delta_{\Phi}=
\{c=d^k\}\cup\Delta_1\cup\Delta_2\cup\Delta_3$, for any $k\geq 3$,  is:

\vspace*{4cm}

\hspace{-1cm}
\fig{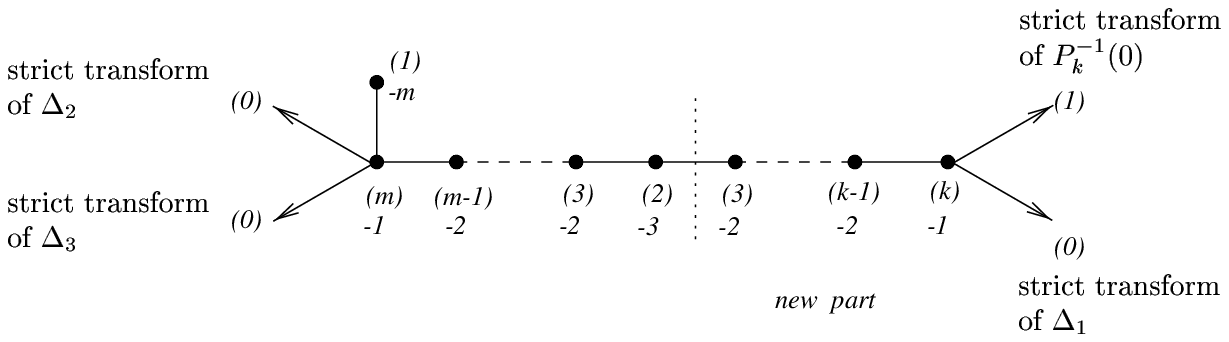}

The reader is invited to construct the corresponding 
diagrams in the case $l=1$. 
In this second case $\Delta_\Phi $ has two components:
$\Delta_1=\{c=0\}$ and $\Delta_2=\{c^{2m-1}=K'd^{2m+1}\}$ (for some 
constant $K'\neq 0$).  Note that for 
 both  $l=1$ and $l=2$ one has $d_{1}=1$. 

\bekezdes{Remark.}{invs} 
Notice that all the germs $f+g^k$ ($k\geq k_0$) define {\em isolated}
 hypersurface singularities.  The name of the series
originates from the  work of I. N.  Iomdin \cite{Io}, 
who studied the case when $g$ was a generic linear form (called 
the ``classical'' case). 

There is a rich literature dealing with the invariants of the germs 
of the (classical or generalized)  Iomdin series. 
Some results compare the invariants of $f+g^k$ with some 
 invariants of $\Phi$. 
For example, using the polar curve method, one can determine the difference 
of the Euler-characteristic of the Milnor fibers
of $f+g^k$ and of the ICIS $\Phi$ using  the multiplicity of 
$\Delta_\Phi$ (for details, see some of the papers
of L\^e  D\~ung Tr\'ang and B. Teissier, e.g. \cite{T1}).

Other results compare some invariant of  $f+g^k$ with the  corresponding 
invariant of the non-isolated singularity $f$. 
Correction term formulae $i(f+g^k)-i(f)$
were obtained for many ``additive''  invariants $i$:
in the classical case by  Iomdin for the case when $i$ denotes  the 
Euler characteristic of the corresponding Milnor fibers
\cite{Io}, by D. Siersma when   $i$ denotes the zeta function \cite{Siersma},
by M. Saito when   $i$ is the spectrum  \cite{Saito},
 and by A. N\'emethi when $i$ is the signature \cite{Ne126,Neta}.
In order to explain why it is more convenient to compute the
{\em correction term} $i(f+g^k)-i(f)$, rather than the 
individual terms $i(f+g^k)$ and $i(f)$, we introduce the following
object. Fix a good representative  of $\Phi$, and 
assume that $\Delta_1=\{c=0\}$.

\bekezdes{Definition.}{wedgedef}
{\em For any integer $M>0$, we define the 
``wedge set'' of  $\Delta_1\subset D^2_\eta$  by
$$\Wedge=\{ (c,d)\in D^2_\eta\  |\  0<|c|<|d|^M \}.$$}
Now, it is easy to verify that 
for a  sufficiently large integer $M>>0$  
 (and with $\eta$ shrinked, if necessary) 
$\Wedge \subset\ D^2_\eta\setminus\Delta_\Phi$. 
In particular, $\Phi$ is a locally trivial fibration over $\Wedge$. 

The point is that in order to determine the individual terms 
$i(f+g^k)$, one needs,  in general, to  understand the usually rather 
complicated monodromy representation associated with $\Phi$. 
On the other hand, 
for some ``nice'' (in some sense ``additive'') invariants (such as  the 
invariants listed above), the correction term
$i(f+g^k)-i(f)$  depends {\em only} on the behavior of $\Phi$
above $\Delta_1$  and  on the representation associated with the fibration
above a wedge set. But this representation is much simpler than the original
monodromy representation, since $\pi_1(\Wedge)={\bfz}^2$. 
In fact, this is exactly the main idea in the construction of the series:
the geometry of the limit object is changed in a small tubular neighborhood 
of a knot. Here the knot is $\Delta_1\cap\partial D^2_{\eta}$, and its 
tubular neighbourhood is  $\overline{\Wedge}\cap\partial D^2_{\eta}$.

Recall that the aim of this work is to compute another
invariant: the resolution graph
of the hypersurface singularities $(X_k,0)=(\{f+g^k=0\},0)$. 
First notice that this invariant is not (really) additive, which makes its
study more difficult. Moreover, our aim  is to obtain
the ``whole'' invariant, not ``just'' a  correction term.
Nevertheless, the wedge set will still  play a
crucial role:  the resolution graphs of all the germs $(X_k,0)$ 
depend only on the  behavior of $\Phi$ over the 
closure of a wedge set.
Moreover, the construction will provide a  ``correction term'' as well. 

Our method is completely different (although obviously not 
absolutely independent)  from the polar curve method, or from 
the usual constructions applied in the cases $i(f+g^k)-i(f)$.
The discussion of some 
connections is the subject of another paper \cite{NSz}, where e.g. 
the relationship between the resolution 
graphs of the singularities $(X_k,0)$ and the monodromy representation
over $\Wedge$ will be developed. 

\section{A graph associated with the  ICIS $\Phi$.}\label{sec5}

\subsection*{A special curve arrangement.}

\vspace{2mm}

\bekezdes{ -- The embedded resolution of $f^{-1}(0)\cup g^{-1}(0)\subset\bfc^3$
 and its stratification.}{ers}
Consider an ICIS $\Phi=(f,g):(\bfc^3,0)\to (\bfc^2, 0)$ 
and  let $(D,0)\subset(\bfc^3,0)$ denote the local divisor 
$(D,0):=((fg)^{-1}(0),0)$. 

Fix an embedded resolution  $r:V^{emb}\rightarrow U$ of 
$(D,0)\subset(\bfc^3,0)$. This means that $U$ is a small representative 
of $(\bfc^3,0)$, $V^{emb}$ is a smooth analytic 3--manifold and 
$r$ is a proper map such that:

(a)\ \
$r|_{{V^{emb}\setminus r^{-1}(Sing(f^{-1}(0))\cup Sing(g^{-1}(0))}}$
is a biholomorphic isomorphism onto its image;

\vspace{1mm}

(b)\ \ the total transform $\de:=r^{-1}(D)$ is a  normal crossing divisor.

\vspace{1mm}

Recall that by the general theory, $r$ is isomorphic over the complement of 
$Sing(D)$. In the above definition, 
 $Sing(f^{-1}(0))\cup Sing(g^{-1}(0))$ is a smaller set than 
 $Sing(D)$, but the intersection
$f^{-1}(0)\cap g^{-1}(0)$ is already transversal off the origin, hence
the above assumption can always be satisfied for a convenient resolution. 

Moreover, we can assume that $U=\Phi^{-1}(D_{\eta}^2)\cap B^3_{\epsilon}$
as well (cf. 4\ \S, I.). 

From (b), for any $p\in \de$, there is a coordinate neighbourhood 
                   $U_p$ of $p$ with local
                   coordinates  $ (u,v,w)$ such that
                   $(fg\circ r)|_{U_p}(u,v,w)=
                   u^{l}\cdot v^{m}\cdot w^{n},$
                   where $l, m, n \geq 0$.

Corresponding to these local equations, $\de$ has a natural stratification:

\bekezdes{Definitions.}{strati} {\em 
To each $p\in\de$ assign the integer $k\in\{1,2,3\}$, if in the above 
local equation there are exactly $k$ non-zero integers among $ l,m,n$.
Denote the stratum corresponding to the number $k$ by $\de^{(k)}$.

Write $\de$ as the union of the following three  subsets:
$\dec =\overline{{ r^{-1}(\{f=0, g\neq 0\})}}$, 
$\ded =\overline{{ r^{-1}(\{g=0, f\neq 0\})}}$ and
$\deo ={ r^{-1}(\{g=0, f= 0\})}$.

A  point $p\in\dee$ is said to be
of type $*$ (where $*=c $ or $d$ or $0$), if $p\in\de_*$.
Similarly,  an irreducible component $B$ of the total transform $\de$
is of type $*$ if its generic point is of type $*$.

A  point $p\in\dek$ is said to be
of type $*_1-*_2$, if $p$ is an intersection point of two local divisors
of type $*_1$ and $*_2$. An  irreducible curve 
$C\subset \overline{ \dek}$ is of type 
$*_1-*_2$ if its generic point is of that type. 
There is a similar definition for the points of the strata $\deh$.}

\vspace{2mm}

In the sequel $\Phi $ denotes a fixed good representative of the ICIS, $r$ is 
an embedded resolution 
as above with  the compatibility  $U=\Phi^{-1}(D_{\eta}^2)\cap B^3_{\epsilon}$.
Define the composition $\wP=\Phi\circ r$. 
Then the restriction 
$\,r^{-1}(\Phi^{-1}(D_\eta^2-\Delta_\Phi)\cap B_\epsilon^3)
\longrightarrow   D_\eta^2-\Delta_\Phi$ of $\wP$ is a locally trivial
           $C^\infty$-fibration, which is 
obviously equivalent with the fibration
induced by $\Phi$ over $D_\eta^2-\Delta_\Phi$. 
Its fibre $\wfcd=\wP^{-1}(c,d)$ is 
called ``a lifted Milnor fibre''.

\bekezdes{The curve $\C$.}{cdef}
Notice that in general (and also in most of the 
concrete examples) it is rather difficult
to find a resolution $r$. Moreover, the codification of its exceptional
divisors and  the corresponding 
intersections and normal bundles can be a rather
difficult problem. Nevertheless,  what we need from this resolution is 
only a 1--dimensional object, the special curve arrangement:

\centerline{$\C=(\dec\cap\deo)\cup(\dec\cap\ded)$.}

\vspace{1mm}

In order to understand the special role $\C$ plays in the geometry of
the series, notice the 
following (see also the comments in \ref{invs}--\ref{wedgedef}).
Since $\{f+g^k=0\}=\Phi^{-1}(\{c+d^k=0\})$ and $\{c+d^k=0\}\subset \overline{
\Wedge}$ (for all $k>M$), all the information that  we need is in the space
above $\overline{\Wedge}$. But $\wP^{-1}(\overline{\Wedge})$ contains 
the 2--dimensional exceptional divisor $\deo=\wP^{-1}(0)$ as well
(a fact which is not very encouraging). However, as it turns out, the relevant 
information is already 
contained in a smaller set,  the closure  of $\wP^{-1}
(\Wedge)$. The intersection of this set with $\deo$ is exactly $\C$
(see below). This shows that from the resolution $r$ we only need a small
tubular neighbourhood of $\C$.

\bekezdes{}{tubn} Let $\C=C_1\cup\cdots\cup C_l$ be the
decomposition of $\C$ into  irreducible components. By construction, each
$C_i$ is either a  smooth curve or it has some 
 self--intersection (double point) singularities.  
Let $T(C_i)$ denote a fixed small 
 open tubular neighbourhood of $C_i$ in $V^{emb}$ ($1\leq i\leq l$), and 
write $T(\C)=\cup_{i=1}^lT(C_i)$. 

\bekezdes{Theorem. -- First characterization of $\C$.}{wedgelift}{\em 
For any open (tubular) neighbourhood $T(\C )$,
there exist a sufficiently large integer $M$  and a  
sufficiently small $\eta$ such that 
$$\wP^{-1}(\Wedge)\ \subset\  T(\C).$$
This follows from the identity:}
\ \ $\overline{\wP^{-1}(\Wedge)}\cap\wP^{-1}(0)=\C.\hfill$

\noi {\em Proof.}\ 
The identity  follows by local verifications. As an exemplification, we
show that if $p\in\dee$ is of type $0$, then 
$p\not\in\overline{\wP^{-1}(\Wedge)}$. Indeed, 
in a local coordinate neighbourhood $U_p$ of $p$, 
$\wP|_{U_p}=(u^m,u^{\mu})$ for some integers  ${\mu},m>0$.
Assume that  $p\in\overline{\wP^{-1}(\Wedge)}$.
Then there exists a sequence $\{z_j\}_{j=1}^\infty$ in $\wP^{-1}(\Wedge)
\cap U_p$ for which $\lim_{j\to\infty}z_j=p$. Hence, if $(u_j,v_j,w_j)$
are the local coordinates of $z_j$, then 
$(u_j,v_j,w_j)\rightarrow 0$.
Moreover, since $z_j\in{\wP^{-1}(\Wedge)}$ it follows that
 $\wP(z_j)\in{\Wedge}$, therefore
$|u_j|^m<|u_j|^{M{\mu}}$.
 This contradicts $M{\mu}-m>0$, for $M>0$ sufficiently large.
All the other local verifications are left to the reader.

The inclusion follows from the identity and the properness of $\wP$.\hfill
$\diamondsuit$

\bekezdes{Corollary.}{connectedhez}  {\em 

a)\  For any open (tubular) neighbourhood $T(\C )\subset V^{emb}$,
there exist a sufficiently large $M$  and a  
sufficiently small $\eta$ such that 
for any $(c,d)\in\Wedge$ the ``lifted Milnor fibre''
 $\wfcd$ is in $T(\C)$.

b)\ For any $p\in\C$ and local neighbourhood $U_p$ of $p$, 
there exist a sufficiently large $M$  and a  
sufficiently small $\eta$ such that 
for any $(c,d)\in\Wedge$ one has $U_p\cap\wfcd\neq\emptyset$.}

\bekezdes{Corollary.}{Cosszefuggo}
{\em The curve  $\C$ is connected.}

\noi{\em Proof.} \ Use (\ref{connectedhez})
and the fact that $\fcd$ (or equivalently, $\wfcd$) is connected. 
\hfill $\diamondsuit$

\bekezdes{Theorem. -- Second characterization of $\C$.}{dlimit} 
{\em Consider the punctured disc 
${\cald}_\eta=\{c=0,d\neq 0\}\cap D_\eta^2\,\subset\,\bfc^2$.
Then}
$$\overline{\wP^{-1}({\cald}_\eta)}\cap \wP^{-1}(0)=\C.$$
Its  proof is similar to that of
(\ref{wedgelift}). The analog of (\ref{connectedhez}) is:
\bekezdes{Corollary.}{dtube} {\em 

a)\  For any  open tubular neighbourhood $T(\C )\subset V^{emb}$,
if $|d|$ is sufficiently small then    \ $\wP^{-1}((0,d))\subset T(\C)$.

b)\ For  any $p\in\C$ and neighbourhood $U_p$ of $p$, if $|d|$ is 
sufficiently small then $U_p\cap\wP^{-1}((0,d))\neq\emptyset$. }

\subsection*{The weighted dual graph of the curve arrangement $\C$.}

\vspace{2mm}

\bekezdes{}{w} 
For each irreducible component $C_i$ of $\C$,  there are exactly two
irreducible components $B_1$ and $B_2$ of the total transform $\de$ for which 
$C_i$ is a component of $B_1\cap B_2$. By the definition of $\C$,
we can assume that $B_1$ is of type $c$ and  $B_2$ is of
type $0$ or $d$. Let $m_{f,B_i}$ (respectively $m_{g,B_i}$) be the 
vanishing order (or multiplicity)
 of $f\circ r$ (respectively of $g\circ r$) along $B_i$
($i=1,2$). Then $m_{g,B_1}=0$, 
$m_{f,B_1}>0$, $m_{g,B_2}>0$, and  $m_{f,B_2}\geq 0$.
Notice that from the  ordered triple $(m_{f,B_1};m_{f,B_2},m_{g,B_2})$
one can recover the local equations for 
$f\circ r|_{U_p}$ and $g\circ r|_{U_p}$ 
in a small coordinate neighbourhood $U_p$ of any point $p\in C_i\cap\dek$.

The components $C_i$ of $\C$  are either  compact (projective) or non-compact.
The compact components are exactly those which are contained in $r^{-1}(0)$.

The non-compact components  form the strict transform of $\{ f =
g = 0\}$. In particular, they are of type $c-d$.
Moreover, since $(f,g)$ is an ICIS, 
the multiplicity of $f$ (resp. of $g$) along the strict transform of
$\{f=0\}$ (resp. of $\{g=0\}$) is one. 
Therefore, the ordered triple $(m_{f,B_1};m_{f,B_2},m_{g,B_2})$
assigned to such an irreducible curve is $(1;0,1)$.

\vspace{2mm}

Now we are ready to define the weighted dual graph 
$\gc$ of $\,\C$. Similarly to the embedded resolution graphs of germs
defined on normal surface singularities
 (cf. \ref{1.1}--\ref{1.4}), 
 the graph  $\gc$ has vertices $\calv$ and edges $\cale$. 

\bekezdes{Definition.}{gc}
The set of vertices $\calv$ consists of 
two disjoint subsets  $\calw$ and $\cala$, where $\calw$ denotes the
non-arrowhead vertices, and $\cala$ denotes the arrowhead vertices.  Let the
non-arrowhead vertices correspond to the compact irreducible curves, 
and the arrowhead vertices to the non-compact irreducible curves of $\C$. 

For any two curves $C_i, C_j\subset\C$,  which
correspond to vertices $v_i, v_j\in\calv$, if $C_i$ and $C_j$
intersect in $l$ points, then let the vertices
$v_i$ and $v_j$ be connected by  $l$ edges
$e_k=\{ (v_i, v_j)_k\}_{k=1,...,l}$ in $\gc$.
Moreover, if a compact component 
$C_i\subset\C$,  corresponding  to a vertex $v_i\in\calw$,
intersects itself, then let each self-intersection point determine 
 a loop (as a special edge)  $(v_i,v_i)$ in the graph $\G_{\C}$.
The edges are not directed, i.e. $ (v_i,v_j)_k=(v_j,v_i)_k$.

\vspace{1mm}
 
{\em Let  the graph $\gc$ be  decorated as follows:}  

$\bullet$\ Let each non-arrowhead
vertex $v_i\in\calw$ have two weights assigned to it:
the ordered triple of integers $(m_{f,B_1};m_{f,B_2},m_{g,B_2})$
   assigned to the {\ic} $C_i$ corresponding to the vertex $v_i$, 
and the genus $g_{i}$ of (the normalization of) $C_i$.

$\bullet$\ 
Let each arrowhead vertex have a single weight, the ordered triple $(1;0,1)$.

$\bullet$\ Let each  edge have a weight $\in\{1,2\}$  determined as follows. 
By definition, any edge $\ e$ corresponds to an intersection point 
$p\in \deh$, hence it is the intersection point of three
local irreducible components of $\de$. The corresponding types of these
components are  $c-*-**$ or $c-c-*$, with  $*\not=c$. 
In the first case let the weight on the edge $e$ be $1$, 
in the second case $2$.

\bekezdes{Summary of notation for $\gc$, and local equations.}{summary}\

\vspace{2mm}

\noindent {\bf Vertices.} \

\noindent \underline{Case 1.}\ \
For all $ w\in\calw$,
\begin{tabular}[b]{p{0.4in}c}
{\includegraphics{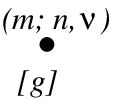}}& means that the
corresponding irreducible 
\end{tabular} curve $C_{w}\subset\C$ is compact and its genus is $g$. 
When $C_w$ is rational,  $[g]$ will  be omitted. 

Furthermore, there is a local coordinate neighbourhood $U_p$ of any
point $p\in C_w\cap\dek$ such that the local irreducible components
$U_p\cap\de= B_1^l\cup B_2^l$ of $\de$ have representation $B_1^l=\{ u=0\}$, 
    $ B_2^l=\{ v=0\}$, and
$$f\circ r|_{U_p}=u^{m}v^{n}, \hspace{7mm} g\circ r|_{U_p}=v^{\nu}\ \ \
\mbox{with } m,\nu>0;\ n\ge 0.$$
Obviously, if  $C_w$ is of type $c-d$ then  $n=0$, and
if $C_w$  is of type $c-0$ then  $n>0$.

\newpage

\noi\alahuz{Case 2.} \begin{tabular}[b]{p{1in}c}
{\includegraphics{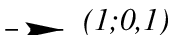}}& denotes an arrowhead $v\in\cala$.
\end{tabular}\newline
 The curve $C_v$ corresponding to $v$ is non-compact.
In a neighbourhood
$U_p$ of a generic point 
$p\in C_v\cap\dek$ we have 
$U_p\cap\de= B_1^l\cup B_2^l$ with $B_1^l=\{ u=0\}$, 
    $ B_2^l=\{ v=0\}$ and
$$ f\circ r|_{U_p}=u, \hspace{1cm} g\circ r|_{U_p}=v.$$

\noindent {\bf Edges.}
An edge $e$  corresponds to an intersection point
$p\in C_{i}\cap C_{j}$ (or a self-intersection of $C_{i}$ if
 $i=j$). In a local coordinate neighbourhood of $p$,
 the irreducible components 
of $U_p\cap\de= B_1^l\cup B_2^l\cup B_3^l$ are given by
$B_1^l=\{ u=0 \},\ B_2^l=\{v=0\},\ B_3^l=\{ w=0\}$.

\vspace{2mm}
\noi\alahuz{Case 1.}

\bs\noi
\begin{minipage}[c]{2.5in}

\vspace{0.4cm}
\begin{tabular}[b]{p{2.5in}}
{\includegraphics{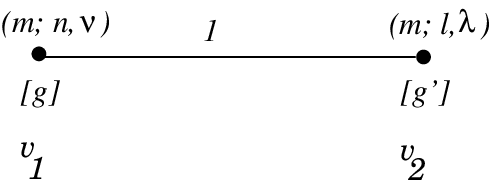}} 
\end{tabular}
\end{minipage}
\begin{minipage}[b]{2.5in}
means  that the unique component of type $c$ is $B_1^l$, and 
$$ f\circ r|_{U_p}= u^{m}v^{n}w^{l}, \hspace{1cm}  
g\circ r|_{U_p}=v^{\nu}w^{\lambda}.$$
\end{minipage}

\vspace{0.2cm}\noi
In fact,  $n=l=0$ means a $(c-d-d)$--\, ;  $n> 0,
l=0$ a $(c-0-d)$--\, ; $n=0$, $l> 0$ a  $(c-d-0)$--\, ; 
and $n> 0, l> 0$ a  $(c-0-0)$--type intersection.



\noi
\begin{minipage}[c]{2.5in}
 If  $v_2$ is replaced by an  arrowhead,\\
 then the corresponding edge is\\
$\ $
\end{minipage}
\begin{minipage}[b]{2.5in}

\vspace{0.9cm}
\begin{tabular}[b]{p{2.5in}}
{\includegraphics{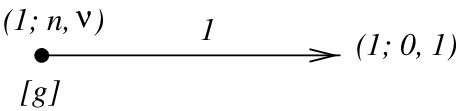}} 
\end{tabular}
\end{minipage}

\noi
 The local equations $f\circ r|_{U_p}$ 
and $g\circ r|_{U_p}$ are as above with $m=\lambda=1$ and  $l=0$.


\sm\noi \alahuz{Case 2.}

\sm\noi
\begin{minipage}[c]{2.5in}

\vspace{0.7cm}\noi
\begin{tabular}[b]{p{2.5in}}
\raisebox{-1ex}{\includegraphics{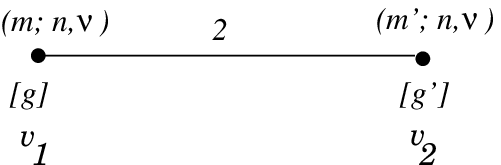}} 
\end{tabular}

\end{minipage}
\begin{minipage}[b]{2.5in}
means that there are two components of type $c$, namely $B_1^l$ and $B_2^l$,
and 
$$
f\circ r|_{U_p}= u^{m}v^{m'}w^{n},\qquad  g\circ r|_{U_p}=w^{\nu}.$$

\end{minipage}

\vspace{0.2cm}\noi
In fact, 
 $l>0$ gives a  $(c-c-0)$--type  intersection, while $l=0$ a $(c-c-d)$--type.


\sm\noi
\begin{minipage}[c]{2.5in}
If $v_2$ is replaced by an arrowhead,\\
 then the corresponding edge is\\
$\ $
\end{minipage}
\begin{minipage}[b]{2.5in}

\vspace{0.7cm}
\begin{tabular}[b]{p{2.5in}}
{\includegraphics{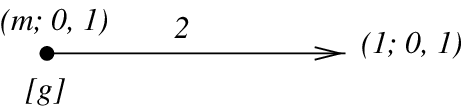}} 
\end{tabular}
\end{minipage}

\noi
 The local equations  $f\circ r|_{U_p}$ and $g\circ r|_{U_p}$ 
are as above with $m'=\nu=1$ and $n=0$.

\bekezdes{A compatibility property of the decorations.}{graphIII} \

\sm\noi
\begin{minipage}[b]{2.5in}
Consider the edge:

\vspace{1.3cm}
\begin{tabular}{p{2.5in}}
\raisebox{-1ex}{\includegraphics{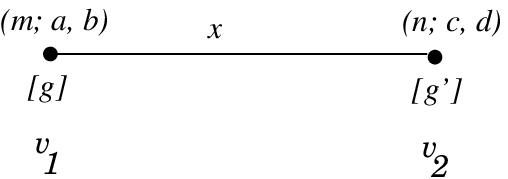}}
\end{tabular}
\end{minipage}
\begin{minipage}[b]{2.5in}
 a.) if $m\neq n$, then $(a,
b)=(c,d)$ and $x=2$;\\
 b.) if $(a,b)\neq (c,d)$, then $m=n$ and $x=1$.

\end{minipage}

\sm\noi In particular, in the cases (a)-(b) above,
the weight of the  edge  is determined by the weights of
the vertices. In all these cases we omit the corresponding edge--decoration.

\newpage
\subsection*{Examples.}

\vspace{2mm}

\bekezdes{Preliminary remarks.}{pr} \

{\bf 1.}\ We emphasize again that in general it is rather long and difficult
to find a resolution $r$. In the literature there are very
explicit resolution algorithms, but they are rather 
involved, and in general very 
slow. Therefore, in almost all of our examples, we used a sequence of some
ad hoc blow ups (following the naive principle: ``blow up the worst
singular locus'') with the hope to obtain a more or less small configuration.
The computations are extremely long, and are not given here. 
(They were, in fact, done with the help of {\em Mathematica}.)
Hence, we admit that for the reader the verification of some of the 
examples listed in the body of the paper
can be a really difficult job. 

\vspace{1mm}

{\bf 2.}\ Since the resolution $r$ is not unique, different resolutions $r$ 
might result in  different graphs $\Gamma_\C$. Moreover, contrary to the
case of resolution graphs of normal surface singularities, in the case of 
$\Gamma_\C$ we cannot claim (at this moment) the existence of a unique
minimal graph. On the other hand, the properties of the graph $\Gamma_\C$ 
(see the 
 next subsection) suggest  that $\Gamma_\C$ inherits a lot of properties
 of  the resolution graphs of normal surface singularities. In particular, 
we expect the existence of some kind of calculus of graphs, which would
provide a method to identify all the possible graphs associated with the
same ICIS $\Phi$. 

\bekezdes{}{hurkosn} Let $f(x,y,z)=x^n+y^n+xyz^{n-2}$ and $g(x,y,z)=z$.
A possible $\gc$ is:

\vspace{8mm}

\vspace{1cm}\hspace{2cm}
\fig{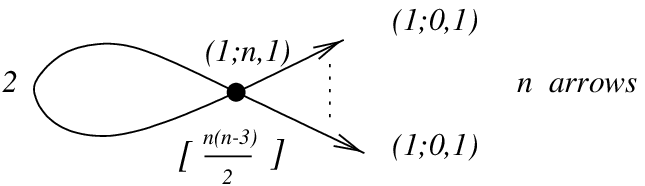}

\vspace{3mm}

More generally, 
if $f$  is given by a homogeneous polynomial, and $g$ is a generic
linear function, 
then the construction of $\gc$ is not difficult. Indeed, first 
blow up the origin. Then consider the intersection curve $C$ of the 
strict transform $St$ of $\{f=0\}$ with the exceptional divisor $E$. 
In the neighbourhood of a  singular point $p$ of $C$, 
 $(St,p)$ has a natural product structure $(C,p)\times (D,0)$, where
$(D,0)$ denotes the disc. Then any embedded resolution of $(C,p)\subset (E,p)$ 
provides an
embedded resolution of $(St,p)$ if the quadratic transformations 
of the plane are replaced by 
the corresponding blow ups along 1--dimensional axes in the 
direction $(D,0)$. 

\bekezdes{}{kistrukkos}
  Let $f(x,y,z)=x^2+y^2$ and $g(x,y,z)=x^n+y^n+z^n$. A possible $\gc$ is:

\vspace{6mm}

\vspace{1cm}\noindent
\fig{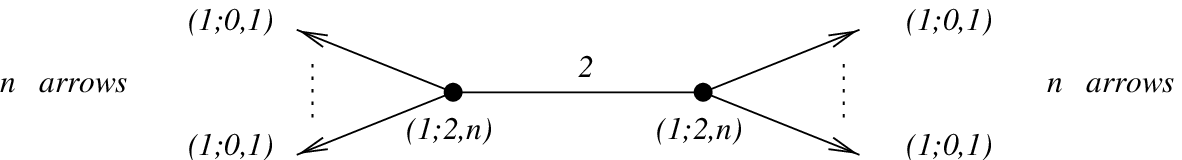}
 
\newpage

\noindent The graph $\G_\C$ is not unique, but depends on the
embedded resolution $r$. Another valid graph $\gc$ 
for the same pair $(f,g)$ is (cf. also with  \ref{extrablowup}):

\vspace{1.6cm}\noindent
\fig{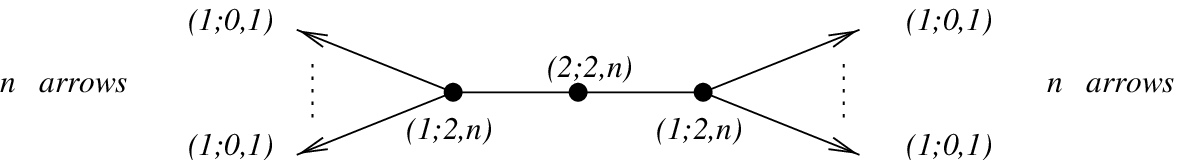} 

\vspace{5mm}

\bekezdes{}{bingo}
  Let $f(x,y,z)=x^2-y^2$ and $g(x,y,z)=y^2+z^3$. A possible $\gc$ is:

\vspace*{2.5cm}

\fig{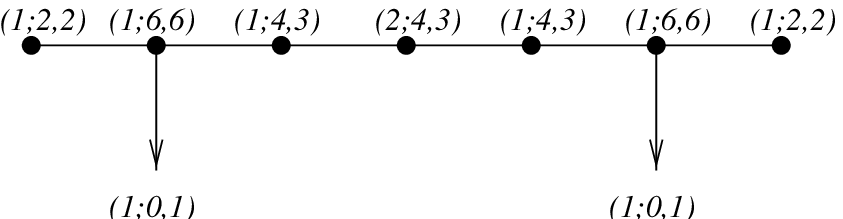}

\vspace{6mm}

\bekezdes{}{andrase}
Consider $f(x,y,z)=x^3+y^2+xyz$ and $g(x,y,z)=z$. A possible $\gc$ is:

\vspace{1.7cm}
 
\vspace{1cm}\hspace{2cm}
\fig{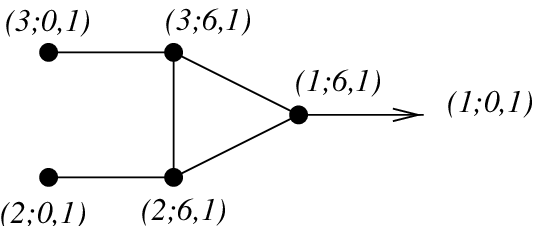}

\vspace{2mm}

\bekezdes{}{duplanyil}
Consider $f(x,y,z)=x^2y^2+z^2(x+y)$ and $g(x,y,z)=x+y+z$. A possible $\gc$ is:

\vspace*{4cm}\hspace{2cm}
\fig{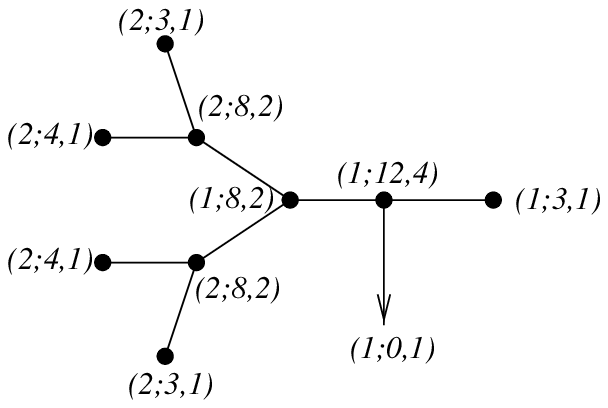}

\bekezdes{}{planecurve}
 Consider $f(x,y,z)=f'(x,y)$ and $g(x,y,z)=z$, where
$f':(\bfc^2,0)\to(\bfc,0)$ is an isolated  plane curve singularity. 
Replacing the quadratic transformations of the plane by blowing ups
along 1--dimensional axis (the local $z$--axis), one obtains the following.
Let $\G(\bfc^2,f')$ denote the minimal embedded resolution graph of the
plane curve singularity $f':(\bfc^2,0)\to(\bfc,0)$.
Then a possible dual graph $\gc$ 
can be obtained from $\G(\bfc^2,f')$ using the following conversion:

$\ $

\vspace{0.7cm}\noi\begin{tabular}{lp{2cm}cp{2cm}}
a non-arrowhead vertex&
\fig{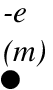}& should be replaced by
&\fig{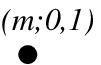}\\
&&\\
an arrowhead vertex &
\fig{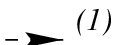}& should be replaced by
&\fig{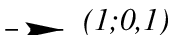}\\
\end{tabular}

\sm\noi
All edges of $\G_\C$ have weight 2 and all vertices $w\in\calw$ of $\G_\C$
have genus $0$.

\vspace*{2mm}

\noindent For example, let $f(x,y,z)=f'(x,y)=x^2-y^3$ and $g(x,y,z)=z$. Then we
have:

\vspace*{1cm}\noi
$\G(\bfc^2,f')$: 

\vspace{1cm} 

\hspace{1.4cm} \fig{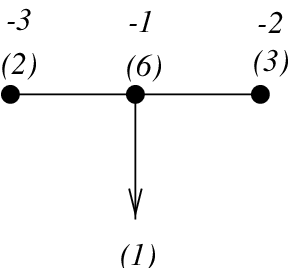}

\vspace{-2cm}

\hspace{6cm} and $\gc$: 

\vspace{1cm}

\hspace{8cm}  \fig{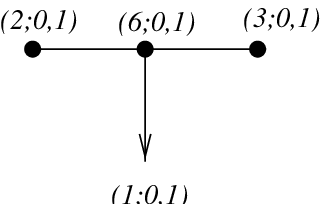}

\subsection*{Some properties of the weighted dual graph $\gc$.}

\vspace{2mm}

\bekezdes{-- The first partition of $\G_\C$.}{gcdecomp}
The vertices of the graph $\gc$ can be divided into two
disjoint sets ${\calv}(\gc)={\calv}^1(\gc)\cup{\calv}^2(\gc)$, where
${\calv}^1(\gc)$ (respectively ${\calv}^2(\gc)$)
 consists  of  the vertices with weight $(m;n,\nu)$ where $m=1$
(respectively $m\geq 2$). 
We will use similar notations for $\calw(\gc)$ and $\cala(\gc)$.

\vspace{1mm}

\noi {\bf The graph $\gce$.} 
Consider the maximal subgraph $\gce$  of $\gc$ which is spanned 
by the vertices $v\in {\calv}^1(\gc)$ and 
has no edges  of  weight $2$ (cf. \ref{extrablowup}). 
Notice that all the arrowheads of $\gc$ are arrowheads of 
$\gce$ as well, i.e. $\calw(\gce)=\calw^1(\gc)$ and $\cala(\gce)=\cala(\gc)$.

If an arrowhead $v$ is supported by an edge  of weight 2, then $v$ becomes an
arrowhead of $\gce$ without any supporting edge. This arrowhead vertex
$v$ (as a completely isolated vertex) forms a connected component of $\gce$. 
For example, in the case of  (\ref{planecurve}), the graph $\gce$
has no edges, and has only one vertex which is an arrowhead.

\vspace{1mm}

\noi {\bf The graph $\gck$.}
The ``complementary'' subgraph $\gck$ is constructed in three steps.
First, consider the maximal subgraph $'\gck$ of $\gc$ spanned 
by the vertices $v\in {\calv}^2(\gc)$. Second, 
in order to obtain $^{''}\gck$, 
 we add some arrowheads (and edges which support these arrowheads) as follows.
Consider an edge $e$ of weight 2 with endpoints $v_2\in \calw^2(\gc)$
and $v_1\in \calv^1(\gc)$. For such an edge, regardless of $v_1$ being an
arrowhead or not, perform the following transformation:

\vspace{2.5cm}

\hspace{-0.5cm}
\fig{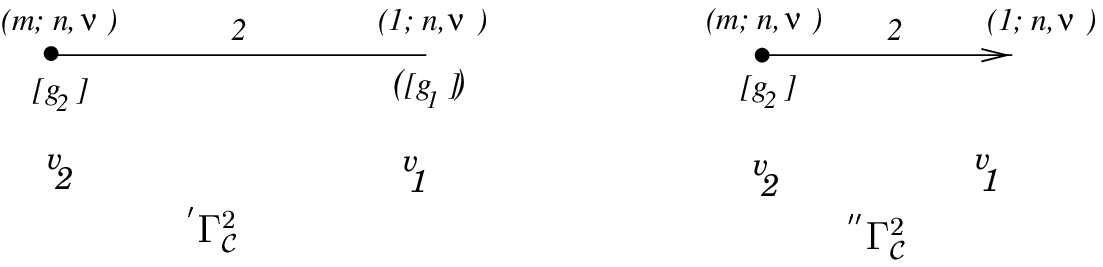}

\noi
Thus $\cale(^{''}\gck)=\cale('\gck)\cup\{\mbox{edges of type $e$}\}$, 
$\calw(^{''}\gck)=\calv^2(\gc)=\calw^2(\gc)$ and $\#\cala(^{''}\gck)=
\#\{\mbox{edges of type $e$}\}$). 

Finally, consider the set of edges in  $\gc$ of weight 2 with both endpoints 
in $\calv^1(\gc)$ (these edges were deleted in the construction of  $\gce$). 
Any edge of this type is transformed into a double arrow:

\vspace*{2mm}
\hspace{5cm}  \fig{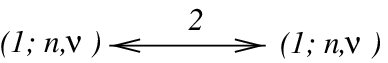}.

\noi
(In general,  this type of graph  codifies the minimal embedded resolution of 
a plane
curve singularity of type $A_1$.) (Cf. also with \ref{extrablowup})

Then, by definition, 
the graph $\gck$ is the union of \,  $^{''}\gck$ and the collection of the
double arrows considered above.  
Obviously, each double arrow forms a connected  component of $\gck$. 

\bekezdes{Example. -- The case of \,2--edges (first case): the double arrows.}{extrablowup} 
Consider an edge of $\gc$ with weight 2 whose  both endpoints are  weighted by
$(1,n,\nu)$ (and arbitrary second weights $[g]$ and $[g']$). 
This edge corresponds
to a point $p\in\deh$ of type $c-c-*$, with $*\in\{0,d\}$.
Moreover, since the first multiplicity of the triplet $(1;n,\nu)$ is 1,
both $c$--type  local irreducible components $B_1^l$ and $B_2^l$
of $\de$ at $p$ are part
of the strict transform of $\{f=0\}$
(cf. \ref{summary}). Their intersection corresponds to the
strict transform of an irreducible component $\Sigma_i$ of $Sing\{f=0\}$ with 
transversal type $A_1$, which has {\em not } been blown
 up during the resolution procedure $r$. 

Performing an additional blow-up along this intersection, we obtain a 
 new embedded resolution $r'$, whose special
curve arrangement $\C'$ will have an additional rational curve.
The dual graph $\gc$ is changed to the
dual graph $\G_{\C'}$ via the transformation:

\vspace{1.5cm}

\fig{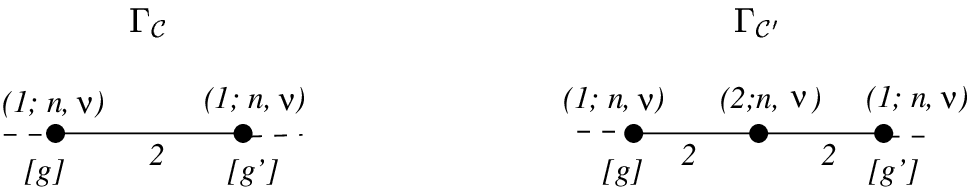}

\sm\noi Therefore, it may be assumed that $\gc$ has no edges of this type. 

\sm
Now, we discuss the properties of $\gc$ and its
subgraphs $\gce$ and $\gck$. 
An immediate consequence of corollary  (\ref{Cosszefuggo}) is the following.

\bekezdes{Proposition.}{} {\em The graph $\gc$  is connected.}

\bekezdes{Remark.}{r2}  Neither $\gce$ nor $\gck$ is connected, in general. 
In fact, the number of connected components of $\gce$ is exactly 
the number of irreducible components of $\{f=0\}$ (cf. \ref{m1}), 
and the number of connected 
components of $\gck$ is exactly the number of irreducible components of
$Sing\{f=0\}$ (cf. \ref{m2tetel}). 
These two facts can easily  be seen on our examples.  E.g.:  in the case of 
(\ref{kistrukkos}), the graph $\gce$ has two components;
and in the case of (\ref{duplanyil})  $\gck$ has two components. 

\bekezdes{-- The garph $\gce$.}{gce} In order to characterize the subgraph
$\gce$ one more notation is needed. Consider the germ 
$(\{f=0\},0)$ and let $n:\{f=0\}^{norm}\to \{f=0\}$ be its normalization. 
Notice that $\{f=0\}^{norm}$ is the disjoint 
 union of some  germs of normal surface
singularities, each germ being the normalization of one of the 
irreducible components of $(\{f=0\},0)$. In particular, 
their number is exactly the number of irreducible components 
of $(\{f=0\},0)$ (see Part I,  \ref{1.7}). 

\bekezdes{Definition.}{m1def} {\em Consider the subgraph $\gce$ of $\gc$.
For any vertex $v\in \calv(\gce)$, replace its weight $(1;n,\nu)$ 
by the weight $(\nu)$. For any non--arrowhead vertex $w\in \calw(\gce)$,
keep its weight $[g]$. Think about this modified
weighted  graph as a weighted  embedded resolution graph of a germ 
of an analytic function defined on
a  normal surface singularity, and about the weights $(\nu)$ as the 
multiplicities of the corresponding function (cf. Part I. \ref{1.2}).
 Calculate ``self-intersection numbers'' $e_w$  ($w\in \calw$)
for each non-arrowhead vertex via the relation (Part I, \ref{1.3}2).
Let the resulting weighted graph be denoted by $\gcce$.

By convention, a dual resolution graph, which consists of only one 
non--arrowhead vertex with weight (1)
 (i.e. it has no other vertices, and no edges),  codifies the dual embedded 
 resolution graph  of a smooth germ defined on a smooth surface, which has 
not been blown up during the resolution.
}

\bekezdes{Theorem.}{m1} {\em 
The  graph  $\gcce$ is a possible  embedded resolution graph of
$$\{f=0\}^{norm}\stackrel{g\circ n}{\longrightarrow} (\bfc,0).$$
In particular, the number of connected components of $\gcce$ 
(which is the same as the  number of connected components of $\gce$) 
coincides with the number of irreducible components of the germ $(\{f=0\},0)$.

Via the above convention, 
a connected component of $\gcce$, which consists of only one non--arrowhead
vertex, corresponds to a  smooth  component of $\{f=0\}^{norm}$ one which the
restriction of  $g\circ n$ defines a smooth map--germ
 (and this smooth component of  $\{f=0\}^{norm}$ 
has not been modified during the resolution).

(Obviously, the subgraph $\gce$ might contain cycles and/or
 non--arrowhead vertices with  $g>0$.)}

\noi {\em Proof.}\ 
Consider the embedded resolution  $r$ and denote by $St$ 
the strict transform of $\{f=0\}$. Then each irreducible component $St_i$
of $St$  is  smooth  and the restriction $r|_{St_i}:St_i\to \{f=0\}_i$ 
is a resolution of the corresponding component $\{f=0\}_i$ of $\{f=0\}$,
which factors through the normalization of $\{f=0\}_i$. Moreover,  
$g\circ r:St_i\to \bfc$ defines a normal crossing divisor on $St_i$. 

Finally, notice that the components $\{St_i\}_i$ are the only $c$--type 
irreducible components of $\de$ on which the vanishing order $m_f$ of 
$f$ is one. The remaining part of the proof consists of verifying
some local
equations and compatibilities of the corresponding graphs, which is left to
the reader. \hfill$\diamondsuit$

\vspace{2mm}

Now we start to analyze the properties of the subgraph $\gck$.

\bekezdes{-- Finer partitions.}{part} Corresponding to the partition 
$\calv^1(\gc)\cup\calv^2(\gc)$ of $\calv(\gc)$,  define  
$\C^k$ as the union of those irreducible components
of $\C$ which correspond to vertices from $\calv^k(\gc)$ ($k=1,2$). 
Obviously, $\C=\C^1\cup\C^2$ and $\C^1\cap \C^2$ is a set of points. 

Now, we define a  finer partition of $\calv^2(\gc)$. 
To any irreducible component  $\Sigma_j$ of $Sing\{f=0\}$, 
\ $1\leq j\leq s$ (cf. \ref{d}),  assign the set 
${\calb}_{\Sigma_j}$ containing those irreducible components $B$ 
of the total transform $\de$ which are of type $c$ and are projected via $r$
{\em onto}  $\Sigma_j$. 
Then define $\csj\subset\C$ as the union of those irreducible components 
$C$ of $\C$ for which $C\subset B$ for some $B\in {\calb}_{\Sigma_j}$. 

\newpage

\bekezdes{Lemma.}{csij} {\em  $\cup_{j=1}^{s}\csj = \C^2$ \ and \ 
$\csi\cap\csj=\emptyset$, $\ $ for any  $i\neq j$, 
\  $\ i,j\in\{1,...,s\}$.}

\noi {\em Proof.} \ The first identity follows from the fact that 
the set of $c$--type irreducible components of $\de$, along which $f$ is
vanishing with order strictly  greater than 1,  
are exactly those components which project 
via $r$ onto some component of $Sing\{f=0\}$. 
For the second part, take a point $p\in\csi\cap\csj$.
Clearly, $p\in\deh$  and it is of type $c-c-*$ with $*\not=c$.
Consider the three local coordinate axes (as in the local description 
\ref{summary}). 
One of them is contained
in a component  $C\subset\csi$, another in a component 
$C'\subset\csj$. Denote the third local coordinate axis by 
$C^*$. Since  it is contained in 
a component $B\in{\calb}_{\Sigma_i}$, it projects via $r$ onto $\Sigma_i$.
But this is true also for $j$, hence $r(C^*)\subset \Sigma_i\cap \Sigma_j=
\{0\}$. But along $C^*$ the function  $g$ is not constant, 
which contradicts $r(C^*)=\{0\}$.  \hfill  $\diamondsuit$

\vspace{2mm}

In the next lemma, $T(\csj)$ denotes the tubular neighbourhood 
$\cup \, T(C_i)$ (cf. \ref{tubn}), where the union is taken 
 over the components $C_i \subset \csj$.

\bekezdes{Lemma.}{ofhoz} {\em 
a)\  For any open (tubular) neighbourhood $T(\csj )\subset V^{emb}$,
there exist a sufficiently small $\gamma>0$ such that for any
point $q\in \Sigma_j-\{0\}\subset \bfc^3$ with $|q|<\gamma$,
$r^{-1}(q) \subset T(\csj)$. 

b)\ For any $p\in\csj$ and local neighbourhood $U_p$ of $p$,
there exist a sufficiently small $\gamma>0$ such that for any
point $q\in \Sigma_j-\{0\}$ with $|q|<\gamma$,
$U_p\cap r^{-1}(q)\neq\emptyset$.}

\noi {\em Proof.} \ The proof is similar to the proof of (\ref{connectedhez})
and (\ref{dtube}) (actually, it can also be deduced from (\ref{dtube})).
\hfill $\diamondsuit$

\bekezdes{Corollary.}{trsconn} {\em 
$\csj$ is connected for all $j\in\{1,...,s\}$. }

\noi{\em Proof.} For any $q\in \Sigma_j-\{0\}$, the fiber $r^{-1}(q)$ 
is connected  by Zariski's Connectivity Theorem (see e.g. \cite{Ha}). 
Then use (\ref{ofhoz}).\hfill$\diamondsuit$

\vspace{2mm}

Therefore, $\C^2$ can be written as a disjoint union of the 
{\em connected} curves $\csj$ ($1\leq j\leq s$). An immediate consequence 
of this is the following.

\bekezdes{Proposition.}{m2tetel} {\em 
There is a one-to-one correspondence 
between the connected components of $\gck$
and the irreducible components of $(Sing\{f=0\},0)$.
The set of non--arrowhead 
vertices in a connected component of $\gck$ corresponds 
exactly to the irreducible components of one of the curves $\csj$.}

\noi ({\em Remark.} Notice that in the case of ``double  arrow components''
of $\gck$, corresponding to the component $\Sigma_j$ (cf. \ref{extrablowup}), 
one has $\csj=\emptyset$.)

\vspace{2mm}

Later  we will show that 
the minimal embedded resolution graph of the transversal singularity
of a connected component $\Sigma_j$
can be reconstructed from the corresponding connected components of $\gck$
(cf. \ref{cts1} and \ref{cts2}). 
But for this we need some further preparation. 

\bekezdes{-- Even finer partitions.}{finpar} Fix a connected component
$\gj$ of $\gck$ corresponding to the irreducible component $\Sigma_j$
$(1\leq j\leq s)$. Its non--arrowhead vertices are denoted by 
$\calw(\gj)$. We define on $\calw(\gj)$ the following equivalence relation.
First, we say that $w_1\sim w_2$ if $w_1$ and $w_2$ 
are connected by an edge of 
weight  1, then we extend $\sim$ to an equivalence relation. The equivalence 
classes $\{K_l\}_l=\{K_l(\gj)\}_l$  determine a partition of $\calw(\gj)$.

Let $\G(K_l)$ be the maximal subgraph of $\gj$ supported by the non--arrowhead
vertices $K_l$. 
For the moment, we keep all the  decorations of the corresponding vertices
of $\G(K_l)$. 
This graph is connected, it has no arrowheads, and its edges are all 
weighted by 1. (If from $\gj$ all the arrowheads and the 
edges of weight 2 are deleted, 
then the disjoint union $\cup_l\, \G(K_l)$ is obtained.)

\bekezdes{Example.}{expar} 
Consider the example given in (\ref{x3y7z4}).
The next figure illustrates the partition of $\G^2_{{\calc},1}$:

\vspace{1.5cm}

\fig{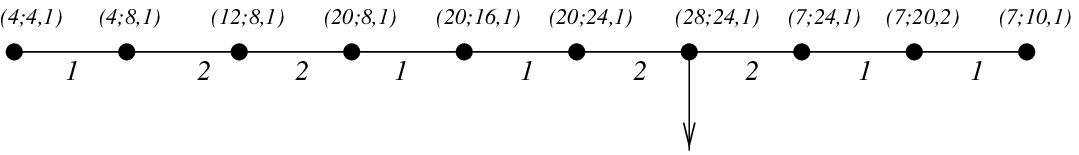}

\vspace{3mm}

\noi The graphs $\G(K_l)$ are:

\vspace{7mm}

\noi\fig{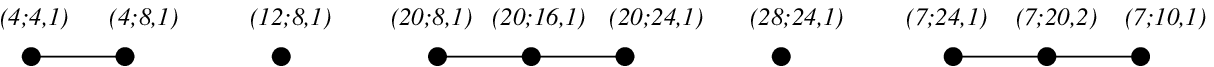}

\bekezdes{}{foly} For a fixed equivalence class
 $K_l=\{w_{i_1},w_{i_2},\ldots, w_{i_t}\}$, consider
the corresponding irreducible curves  $\C_{i_1},\C_{i_2},$ $\ldots,
 \C_{i_t}$ of
$\C$. By construction, there exists an irreducible component $B(K_l)\in
{\calb}_{\Sigma_j}$ which contains all of them. Moreover,
the union $\C(K_l):=\C_{i_1}\cup\cdots\cup\C_{i_t}$ is a connected curve.
Let $T(K_l)$ denote  a small tubular neighbourhood of $\C(K_l)$ in $B(K_l)$. 

Finally, consider the weights $(m_{i_k};n_{i_k},\nu_{i_k})$
 associated with $w_{i_k}$ (or with $\C_{i_k}$) \,  ($1\leq k\leq t$).
Again, by the  definition of the partition, $m_{i_1}=\cdots= m_{i_t}$.
We denote this integer by $m(K_l)$;  it is exactly  the vanishing order of
$f\circ r$ along $B(K_l)$. On the other hand, the local 
equations show that the restriction $g\circ r|_{T(K_l)}:T(K_l)\to \bfc$
provides the principal divisor $(g\circ r)^{-1}(0)=
\sum _k\, \nu_{i_k}\, \C_{i_k}$ in $T(K_l)$. 

Therefore, the divisor (multiple curve)
$\sum _k\, \nu_{i_k}\, \C_{i_k}$ in $T(K_l)$
 can be interpreted as a central fiber
 of  a proper analytic map (here $T(K_l)$ can be changed to 
$(g\circ r|B(K_l))^{-1}$(small disc)).\\ 
It is convenient to introduce the integer $\nu(K_l):=g.c.d.(
\nu_{i_1},\ldots,\nu_{i_k})$. 

\bekezdes{Lemma.}{genfib} {\em The generic fiber of 
$(g\circ r)|T(K_l)$ is a disjoint union of rational curves.}

\noi {\em Proof.} \ Fix $d\not=0$ with $|d|$ sufficiently small. Then 
$g^{-1}(d)$ intersects $\Sigma_j$ in $d_j$ points $\{q_i\}_i$ (cf. \ref{d}). 
Then $((g\circ r)|T(K_l))^{-1}(d)=\cup_i\,( r^{-1}(q_i)\cap T(K_l))$.
But $r^{-1}(q_i)$ is a fiber of an embedded resolution of the 
transversal singularity associated with $\Sigma_j$, hence each 
irreducible curve  in $r^{-1}(q_i)$ is rational. \hfill $\diamondsuit$

\bekezdes{-- General properties of morphisms whose generic fiber is rational.}{ge}
First we recall the following fact (cf. \cite{GH}, page 554):

\vspace{1mm}

\noi {\em  If $S$ is a minimal smooth surface, $D$ a disc in $\bfc$, 
and $\pi:S\to D$  any proper holomorphic map  whose generic fiber
$\pi^{-1}(d)$ (\,$0\not=d\in D$)  is irreducible and rational, then 
$\pi $ is a (trivial) ${\bf P}^1$--bundle over $D$.}

\vspace{1mm}

Now, assume that $S$ is not minimal, and the generic fiber of
$\pi:S\to D$  is a disjoint union of (say, $N$) rational curves. Then 
by the Stein Factorization theorem (see e.g. \cite{Ha}, page 280),
(and shrinking $D$ if necessary), there exists a
map  $b:D'\to D$ given by  $z\mapsto z^N$, and 
$\pi':S\to D'$ such that $\pi=b\circ \pi'$, and the generic fiber of 
$\pi'$ is irreducible and rational. 
Since the central
fibers of $\pi $ and $\pi'$ are the same, it follows  from the above fact that 
the central fiber of $\pi$ can be blown down 
successively until an irreducible  rational curve is obtained.

This  discussion has the following consequences:

\bekezdes{Proposition. -- Properties of the graph $\G(K_l)$.}{gkl} \ {\em 

a)\ The graph $\G(K_l)$ is a tree with all $g_{w_k}=0$. In particular, 
all the irreducible components of \, $\C^2$ are rational curves.

b)\ From the integers $\{\nu_{i_k}\}_{k=1}^t$ one can deduce the 
self--intersections $C_{i_k}^2$ of the curves $C_{i_k}$ in $B(K_l)$ as follows.
First notice that the intersection matrix $(C_{i_k}\cdot C_{i_{k'}})_{k,k'}$
(where the intersections are considered in $B(K_l)$) is a negative 
semi--definite matrix with rank $t-1$, and the ``central divisor'' 
$\sum_k\nu_{i_k}C_{i_k}$ is one 
element of its kernel (cf. \cite{BPV2}, page 90). The intersections 
$C_{i_k}\cdot C_{i_{k'}}$ for $i_k\not=i_{k'}$ can be read from the graph
$\G(K_l)$ considered as a dual graph. The self--intersections can be
determined from the relations $(\sum_k\nu_{i_k}C_{i_k})\cdot C_{i_{k'}}=0$
(cf. also with Part I, \ref{1.3}2). In particular, if
$t=1$, then  $C_{i_1}^2=0$. If $t\geq 2$, then the graph (curve configuration)
is not minimal; if we blow down succesively all the $(-1)$--curves we obtain 
a rational curve with self--intersection zero. This gives a complete 
classification  of the shapes of the possible graphs $\G(K_l)$ (decorated with
the set of integers $\{\nu_{i_k}\}_k$). 

c)\ The number of irreducible (equivalently, connected) components of the 
generic fiber of $(g\circ r)|T(K_l)$ is $\nu(K_l)=g.c.d.(\nu_{i_1}, \ldots,
\nu_{i_k})$. 

The fact that each irreducible component of the generic fiber is rational
can be translated into the relation:
$$2\cdot \nu(K_l)=\sum_{k=1}^t\,\nu_{i_k}(2-\delta_{w_{i_k}}),$$
where $\delta_{w_{i_k}}$ is the number of  vertices adjacent to $w_{i_k}$ 
in $\G(K_l)$.}

\bekezdes{Example.}{expar2} Consider the example (\ref{expar}) (based on 
\ref{x3y7z4}). Then the partition $\{\G(K_l)\}_l$ of $\G_{{\C},1}^2$ 
decorated with the integers $(\nu_{i_k})_k$ and  the self--intersection numbers
$C_{i_k}^2$ is the following:

\vspace{1cm}

\fig{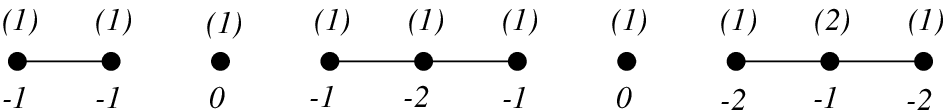}

\bekezdes{}{back} Now we return to the connected component 
$\gj$ of $\gck$ corresponding to $\Sigma_j$ ($1\leq j\leq s$). 
In (\ref{finpar}--\ref{expar2}) we introduced the partition
$\{\G(K_l)\}_l$. The arrowheads of $\gj$ as well as the edges connecting 
different  classes $\{\G(K_l)\}_l$ in $\gck$
 are 2--edges. Their detailed presentation is the subject
 of the next paragraph.

\bekezdes{Example. -- The case of \, 2--edges (second case).}{2edg} \
Consider an edge $e$
of $\gc$ of weight 2 having endpoints  $v_1$ and $v_2$
such that $v_2$ is weighted with the triple $(m';n,\nu)$ with $m'\geq 2$. 
Consider the local coordinates $(u,v,w)$ as in (\ref{summary}),
with $B_1^l=\{u=0\}, \ B_2^l=\{v=0\}$ and $B_3^l=\{w=0\}$. 

First assume that the first entry of the weight
$(m;n,\nu)$ of  $v_1$ also satisfies  $m\geq 2$, i.e.
both $v_1$ and $v_2$ are non--arrowheads of $\gj$. Then $C_{v_1}\cap
U_p$ and $C_{v_2}\cap U_p$  are the intersections $B_1^l\cap B_3^l$
and $B_2^l\cap B_3^l$. Their intersection point $p$ is codified in the edge 
$e$.  From the local equation $g\circ r|U_p=w^{\nu}$, we obtain that 
in the deformation $g\circ r|B_1^l\cup B_2^l\to \bfc $ the intersection
point $p$ splits into $\nu$ points. 

If the first entry $m$ of the weight $(m;n,\nu)$ of $v_1$ is $m=1$, then 
the discussion is similar.

In both cases,
set $C^*_e:=B_1^l\cap B_2^l$. Then $r|C^*_e:C^*_e\to \Sigma_j$ is finite.
Denote its degree by $d(e)$. Obviously $deg(r|C^*_e)\cdot deg(g|\Sigma_j)=
deg(g\circ r|C^*_e)=\nu$. Sometimes it is more convenient (and precise)
to use the notation $\nu(e)$ for $\nu$. With these notations:
$$d(e)\cdot d_j=\nu(e).$$

Now we are able to characterize  the structure of the subgraph
$\gj$. The final goal is to reconstruct the embedded resolution graph
of the transversal singularity associated with the branch $\Sigma_j$
from the weighted graph $\gj$.
First we start with a construction. The reader is invited to
review the theory of cyclic covering of graphs (Section \ref{sec1}, Part I).

\bekezdes{-- A covering construction from  $\gj$.}{covgj} \ 
Fix an irreducible component $\Sigma_j$ of $Sing\{f=0\}$ ($1\leq j\leq s$),
let $T\Sigma_j$ be the (equisingular type of the) transversal singularity
associated with $\Sigma_j$ (cf. \ref{transversdef}). 
Recall that $deg(g|\Sigma_j)=d_j$.

If $(H,q)$ is a
transversal slice as in (\ref{transversdef}), then $r$ above $(H,q)$ 
determines a resolution of the transversal plane curve singularity
$(H\cap \{f=0\},q)\subset (H,q)$. Its weighted 
dual embedded 
resolution graph  (for a  definition, see Part I, \ref{1.2}) is denoted by
$G(T\Sigma_j)$. Since in local coordinates it is easier to work with 
the pullback of $g$, it is convenient to replace the single point 
$q\in \Sigma_j-\{0\}$ by the collection of $d_j$  points 
$g^{-1}(d)\cap \Sigma_j$ (where $|d|$ is small and non-zero). 
The dual weighted graph associated with the curves situated 
above these points consists of  exactly $d_j$ identical
copies of $G(T\Sigma_j)$, and it is  denoted by  $d_j\cdot G(T\Sigma_j)$.

Comparing the curves  $r^{-1}(g^{-1}(d)\cap \Sigma_j)$ and 
$\C_{\Sigma_j}$ by  the corresponding local equations, and using the 
results of (\ref{gkl}), especially part (c), and (\ref{2edg}), we obtain 
a cyclic covering of graphs 
$$p:d_j\cdot G(T\Sigma_j)\to \mbox{\{a base graph\}}_j,$$
where the base graph and the covering data can be determined from $\gj$.
This is given in the next paragraphs.

\vspace{2mm}

\noindent {\bf The base graph} will be denoted by  $\gj/\sim$.
It  is obtained from $\gj$ by collapsing
it along edges of weight 1. 

More precisely, each subgraph $\G(K_l)$ is  replaced by a non--arrowhead 
vertex. If two subgraphs $\G(K_l)$ and $\G(K_{l'})$ are connected by 
$k$ 2--edges  in $\gj$, then  the corresponding
vertices of $\gj/\sim$ are connected by $k$ edges. 
(In fact, from \ref{tree} will follow that $k\leq 1$.) 
If the non--arrowhead vertices of $\G(K_l)$ support $k$ arrowheads altogether, 
then on the corresponding non--arrowhead  vertex of 
$\gj/\sim$ one has exactly $k$ arrowheads. 

Since $\gj$ is connected, it is obvious that $\gj/\sim$ is connected as well.

\vspace{2mm}

\noi {\bf The covering data of}  $p:d_j\cdot G(T\Sigma_j)\to \gj/\sim$.

First we recall that the {\em covering data } of a graph $\G$ is a collection
 of positive integers  $\{n_v\}_{v\in \calv(\G)}$ and 
$\{n_e\}_{e\in \cale(\G)}$, such that for each edge $e\in \cale(\G)$
with endpoints $v_1$ and $v_2$, \ $n_e=d_e\cdot l.c.m(n_{v_1},n_{v_2})$
for some integer $d_e$. For details, see Section \ref{sec1} (Part I), 
or at least (\ref{4.4}).

Now, we define a covering data for $\gj/\sim$. It is provided by the 
third entries $\nu$  of the weights $(m;n,\nu)$ of the vertices of
$\gj$. 

For any non--arrowhead vertex  $w$ of $\gj/\sim$, which corresponds to $K_l$ 
in  the above construction, set $n_w:=\nu(K_l)$.
For any arrowhead vertex $v$ of $\gj/\sim$, which corresponds to
an arrowhead of $\gj$ with weight $(1;n,\nu)$, set $n_v:=\nu$. 
For any edge of $\gj/\sim$, which comes from a 2--edge $e$ of $\gj$
with endpoints with weight $(*;n,\nu)$, set $n_e:=\nu \, (=\nu(e))$. 

For simplicity (and suggestively) the covering data (which in 
Section \ref{sec1} (Part I) 
is denoted by ${\bf (n,d)}$), will here be denoted by 
$\bnu$. 

\vspace{2mm}

The degeneration of $r^{-1}(g^{-1}(d)\cap \Sigma_j)$ into $\C_{\Sigma_j}$
provides the (already announced) result:

\bekezdes{Theorem. -- Characterization of the transversal singularities
(first part).}{cts1} \ {\em 
For any $1\leq j\leq s$, there exists
an embedded resolution graph $G(T\Sigma_j)$  of the transversal
singularity associated with the branch $\Sigma_j$, and 
 a cyclic covering of graphs
$$p:d_j\cdot G(T\Sigma_j)\to \gj/\sim$$
with covering data $\bnu$ (described above) and with the compatibility of the 
arrowheads: $\cala(d_j\cdot G(T\Sigma_j))=p^{-1}(\cala(\gj/\sim))$
(cf. \ref{4.20} (1)). 

In particular, by (Part I, \ref{cgcg}), 
(and from the fact that $G(T\Sigma_j)$ is a 
tree,  being the graph of a plane curve singularity), one obtains that
$\gj/\sim$ is a connected tree.}

\bekezdes{Remark.}{nemtriv} It is not easy to find concrete examples,
when the covering data $\bnu$ is not trivial (i.e. when  the integers 
 $n_v$ and  $n_e$ are not all the same). 
Example (\ref{duplanyil})
 shows such an example (see \ref{nemtrivi} for the detailed
discussion). 

\bekezdes{}{tovabb} The above theorem (\ref{cts1}) has some important 
consequences.  Before we formulate the first, we invite the reader to
verify the following fact.

\vspace{1mm}

\noi {\em If $\G$ is a tree, and 
$p:G\to \G$ is a cyclic covering of \, $\G$ with covering data $\{n_v\}_{v\in 
\calv(\G)}$, then the number of connected components of $G$ is
$g.c.d.\{n_v|v\in \calv(\G)\}$.}

\vspace{1mm}

This implied in (\ref{cts1}), where $G(T\Sigma_j)$ and 
$\gj/\sim$  are connected, gives:

\bekezdes{Corollary.}{dj} {\em \

(a)\ 
$d_j=g.c.d.\{\nu_w|w\in \calw(\gj)$,  where $(m_w;n_w,\nu_w)$ is the weight of
$C_w$\}.

(b)\  Let $\cale_j$ be the set of 2--edges of $\gc$ whose endpoints
are non--arrowhead  vertices such that
 one of the endpoints  is  in $\calw(\gce)$, the other in $\calw(\gj)$. 
(By construction, $\cale_j$ is exactly the index set of the arrowhead
vertices of $\gj$.)
Let $\#(T\Sigma_j)$ be the number of irreducible components of the 
transversal singularity $T\Sigma_j$. Then:
$$d_j\cdot \#(T\Sigma_j)=\sum_{e\in \cale_j}\, \nu(e).$$
}

\vspace*{-2mm}

\bekezdes{Example.}{nuketto} In the example  (\ref{duplanyil}), 
$Sing\{f=0\}$ has two components. For both values $j\in \{1,2\}$, the
transversal type of $\Sigma_j$ is $A_1$, $\#(T\Sigma_j)=2$, $d_j=1$,
$\#\cale_j=1$, and $\nu(e)=2$ for $e\in \cale_j$.

\bekezdes{Example.}{ciklus} If $\#(T\Sigma_j)\not=1$, then even if $\gce$ and
all $\gj$'s are trees, it is possible that $\gc$ has some cycles. 
For example, if $f=x^5+y^2+xyz$ and $g=z$, then $Sing\{f=0\}=\{x=y=0\}$
and has transversal type $A_1$. $\gce$ is a tree (representing a 
surface singularity of type $A_2$), $\gj$ is a tree, but $\gc$ is:

\vspace{1.7cm}

\fig{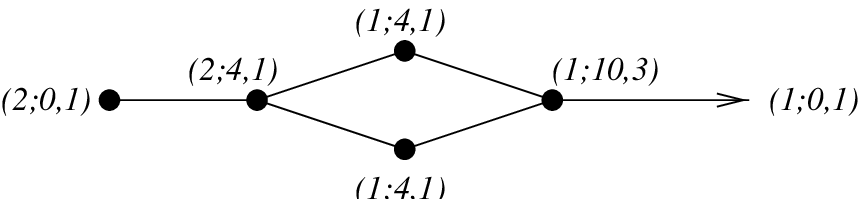}

\bekezdes{Remark.}{stricttra} During the construction of $\gce$ we deleted all 
2--edges (with at least one end--point in $\calw(\gce)$), 
but these  2--edges have a geometric meaning in the case of $\gce$ 
as well. Take such an edge $e$ with endpoints $v_1$ and $v_2$.
If $v_1\in \calw^1(\gc)$ but $v_2\not\in \calw^1(\gc)$, then add to $\gce$ 
an arrowhead supported by  $v_1$. 
If both $v_1$ and $v_2$ are in $\calw^1(\gc)$, then add to $\gce$ 
two  arrowheads, one of them supported by $v_1$, the other by $v_2$.
Let ${\cals}$  denote the collection of these arrowheads. 
Then ${\cals}$ represents (in a similar way as the arrows of
 any dual graph) the strict 
transforms  of the singular locus $Sing\{f=0\}$ in the resolution 
$r|_{\coprod St_i}:\coprod St_i\to \{f=0\}$ (cf. \ref{m1}).
In particular, $\sum_j\#(T\Sigma_j)\geq \#{\cals}$. 

In general the inequality can be strict.  For example, in the case 
(\ref{duplanyil}) (cf. also with \ref{nuketto}), 
$\sum_j\#(T\Sigma_j)=4$ but $\#{\cals}=2$. In this case, the strict transform
$C^*_j$ of $\Sigma_j$ is irreducible, even if the transversal type of
$\Sigma_j$ has two components (but $deg(C^*_j\to \Sigma_j)=2$). 

\vspace{2mm}

\noi Now, we continue the list of consequences of theorem (\ref{cts1}). 
By theorem (\ref{gkl}) $\G(K_l)$ is a tree for each $l$. Since $\gj/\sim$
is a tree as well, one gets:

\bekezdes{Corollary.}{tree} {\em  \
a)\ $\gj$ is a connected tree.

\vspace{1mm}

b)\ 
Moreover, since $\gj/\sim$ is a tree, corresponding to the/any  covering data
$\bnu$, (by Part I, \ref{4.16})
there is \underline{only one} cyclic graph covering 
$p:G\to \gj/\sim$ (up to isomorphism). Therefore, the shape of the graph
$G(T\Sigma_j)$ is completely determined by the
 weighted graph $\gj$.}

\vspace{2mm}

In the sequel we discuss the compatibility of the decorations:
we will show that one can recover all the decorations of $G(T\Sigma_j)$
from the decorations of $\gj$. 

First we recall (cf. Part I, \ref{1.2}) 
that the embedded resolution graph $G(T\Sigma_j)$
has three type of decorations: (a)\ 
each vertex $v\in \calv(G(T\Sigma_j))$ has a {\em multiplicity} $(m_v)$;
(b)\ 
each vertex $w\in \calw(G(T\Sigma_j))$ has a {\em self--intersection} $e_w$;
(c)\ 
each vertex $w\in \calw(G(T\Sigma_j))$ has a {\em genus} $[g_w]$.

\bekezdes{Theorem. -- Characterization of the transversal singularities
(second part).}{cts2} \ {\em 
For any $w\in \calw(G(T\Sigma_j))$ corresponding to $K_l$, $m_w$ is exactly
$m(K_l)$; and 
for any arrowhead $v\in \cala(G(T\Sigma_j))$: \,  $m_v=1$.

By (Part I, \ref{1.3}2), 
each self-intersection number $e_w$ can be determined from the set 
of multiplicities $\{m_v\}_{v\in\calv}$. 
Finally, $g_w=0$ for each $w$. 

\vspace{1mm}

In particular, the weighted dual (embedded) resolution graph $G(T\Sigma_j)$
 of the transversal singularity associated with $\Sigma_j$  can be completely 
determined from the weighted graph $\gj$.}

\bekezdes{Remark.}{conditions} As we already emphasized in (\ref{gkl} b--c),
the collection of integers  $\{\nu_w\} \, (w\in \calw(\gj))$ satisfies
serious compatibility restrictions. 
Moreover, since in the cyclic covering $d_j\cdot G(T\Sigma_j)\to
\gj/\sim$ the covering graph  $d_j\cdot G(T\Sigma_j)$ has no cycles, 
this imposes  some  additional restrictions for the 
integers $\{\nu_w\} \, (w\in \calw(\gj))$. 

\vspace{2mm}

The results of this section  and theorem (Part I, \ref{theo2})
 culminate in the following corollary 
which is  crucial for the algorithm presented in the next section. 

\bekezdes{Proposition.}{unique} {\em Up to isomorphism of cyclic covering of
graphs (with a fixed  covering data), there is only one cyclic covering of the
graph $\gc$ provided that the covering data satisfies $n_v=1$ for any 
$v\in \calv^1(\gc)$.}

\subsection*{ Illustrations of properties of the graph $\gc$.}

\vspace{2mm}

\begin{figure}[h]
\hspace{-11.5cm}\fig{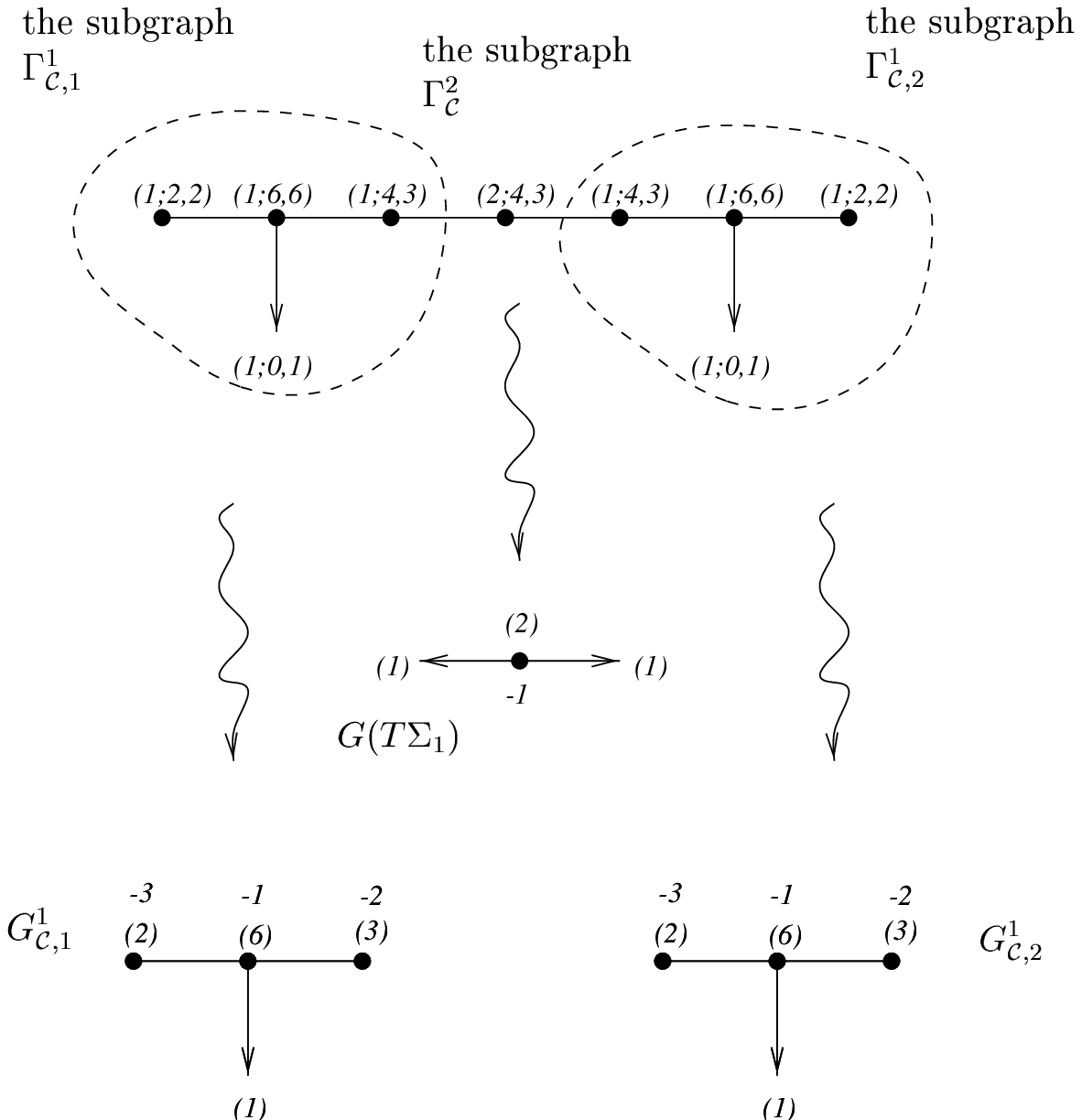}

\vspace{1.5cm}
\caption{}\label{nagybingokep}
\end{figure}

\bekezdes{Example.}{1} 
Consider the ICIS  $(f,g)$, where $f=x^2-y^2$ and $g=y^2+z^3$
(cf. \ref{bingo}). Figure \ref{nagybingokep} 
 shows the decomposition of a possible
$\gc$ into $\gce$ and $\gck$. The graph $\gce$ has two components 
$\G^1_{\C,1}$ and $\G^1_{\C,2}$ corresponding to the two smooth 
components of the (normalization) of $\{f=0\}$. 
The corresponding resolution graphs
$G^1_{\C,1}$ and $G^1_{\C,2}$ provided by 
(\ref{m1def}) represent both an $A_2$ plane curve singularity corresponding 
to the restriction of  $g$ to the irreducible components of $\{f=0\}$. 
The transversal singularity of $Sing\{f=0\}$  is of type $A_1$
with graph $G(T\Sigma_1)$. 

\bekezdes{Example.}{nemtrivi} Consider the example (\ref{duplanyil})
(cf. also \ref{nuketto} and \ref{stricttra}).
 Then $Sing\{f=0\}$ has two components. For both $j\in\{1,2\}$ fixed,  
$d_j=1$ and there is only one equivalence class $K_l$.
Therefore,  $\gj/\sim$ consists of only one non-arrowhead which supports 
exactly one  arrowhead.  But $G(T\Sigma_j)$ has one non--arrowhead which 
supports two  arrowheads  (representing an $A_1$) singularity. Therefore,
in this case the covering $G(T\Sigma_j)\to \gj/\sim$ is not a product
(or identity). 

\begin{figure}[h]
\hspace{-12cm}\fig{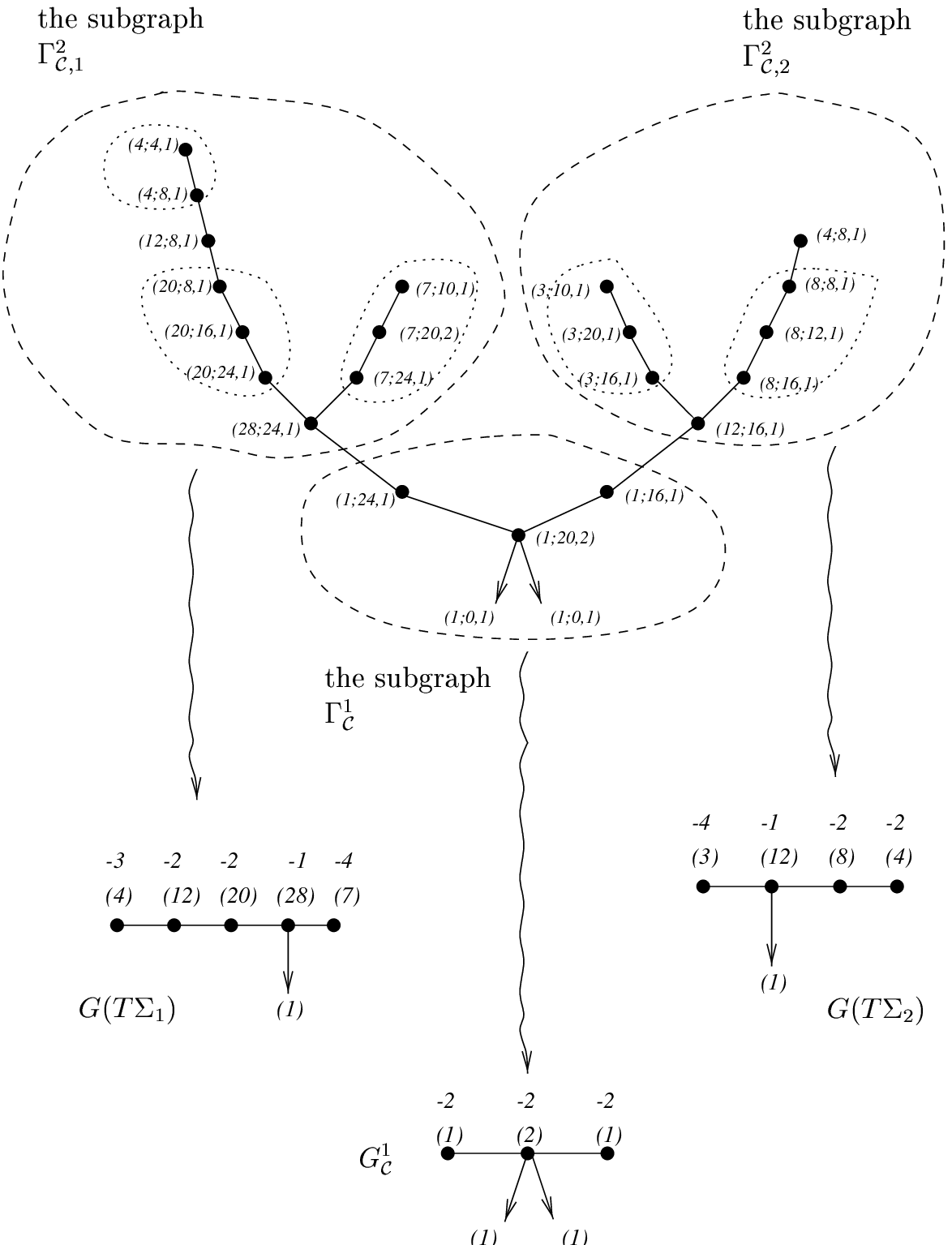}

\vspace{1cm}
\caption{}\label{nagypeldakep}
\end{figure}

\bekezdes{Example.}{x3y7z4}  Consider the ICIS 
 $(f,g)$, where $f=x^3y^7-z^4$ and $g=x+y+z$. In this example,
the reader can verify that the example provides non--trivial 
partitions $\{K_l\}_l$. The case $j=1$ is explained in 
(\ref{expar}) and (\ref{expar2}). Figure \ref{nagypeldakep}
shows the graph $G^1_{\C}$, as well as  $G(T\Sigma_1)$  and
 $G(T\Sigma_2)$ (corresponding 
to the two transversal singularities), each of which is
obtained from a possible $\G_\C$.

\begin{figure}[h]
\hspace{-12cm}\fig{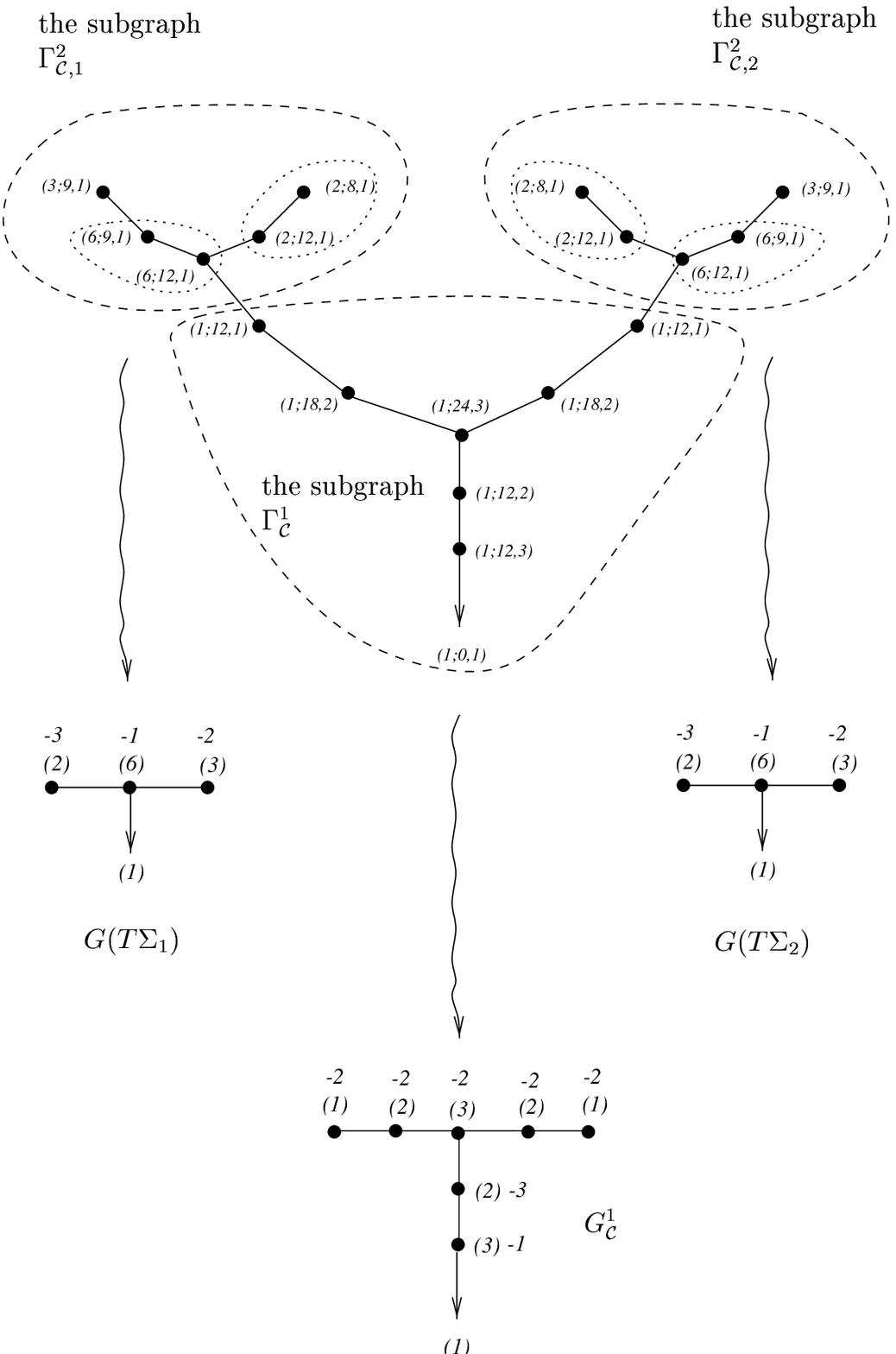}

\vspace{1cm}
\caption{}\label{ujpeldakep}
\end{figure}

\bekezdes{Example.}{ujpelda} Consider the pair $f=y^3+(x^2-z^4)^2$ and $g=z$.
Using the same notations as before,
 Figure \ref{ujpeldakep}
shows the graph $G^1_{\C}$, as well as  $G(T\Sigma_1)$  and
 $G(T\Sigma_2)$ (corresponding 
to the two transversal singularities), each of which is
obtained from a possible $\G_\C$.

\newpage $\ $

\vspace*{13cm}

\newpage $\ $
\vspace*{16cm}

\newpage $\ $

\vspace*{17cm}

\newpage $\ $

\section{The resolution graph  of $\mathbf{\{f+g^k=0\}}$.}\label{sec6}


\vspace{2mm}

Consider a pair of germs 
$f,g:(\bfc^3,0)\rightarrow (\bfc,0)$, such that $f$ has a
1-dimensional singular locus and the pair $\Phi=(f,g)$ forms an
ICIS. Let $k>0$ denote a sufficiently large integer (cf. \ref{goal2}).
Then 
$$(X_k,0):=(\{f+g^k=0\},0)\,\subset\,(\bfc^3,0)$$
\noi is an {\em isolated} hypersurface singularity. 
The algorithm presented in this section provides the minimal
resolution graph $\gxk$ of $(X_k,0)$.

\vspace{2mm}

Actually, our algorithm gives even more: we  construct the embedded
resolution graph $\gxg$ 
of $g|X_k:(X_k,0)\to (\bfc,0)$. In particular, the graph  will have 
three types of decorations: multiplicities, genera and self--intersections,
and also  arrowheads representing the strict transforms of $\{g|X_k=0\}$
decorated with multiplicities as well. 
If we delete from the weighted graph $\gxg$
 the arrowheads and multiplicities
then we obtain the wanted resolution graph $\gxk$ of $(X_k,0)$. 


\bekezdes{-- A short preview of the algorithm.}{algshort} 
The algorithm starts with the 
the construction of the  weighted dual graph $\G_\C$
 associated with an embedded resolution $r$ of
$(D,0):=(\{fg=0\},0)\subset(\bfc^3,0)$ (cf. \ref{ers}).
The graph $\G_\C$ is not unique,
it  depends on the embedded resolution $r$. However, the algorithm works for
{\it any} dual graph $\G_\C$.

Obviously, in order to obtain $\gc$, the optimal is  to have a detailed
description of the resolution $r$. But, let us emphasize again, in order to
identify the necessary
 information for $\gc$ from the resolution, only a
small part of the resolution is needed: some  information  from a small
neighborhood of the  special curve arrangement $\C$.

This construction is  clarified  completely in (\ref{cdef} and  \ref{gc}). 
In the sequel we assume  that a weighted 
graph $\gc$ have  already been   constructed.    
Then our algorithm  provides, in a {\em purely combinatorial way},
the weighted dual embedded resolution graph $\gxg$ 
from the  weighted dual graph $\G_\C$ and the integer $k$.

The reader is invited to review the definitions regarding
cyclic covering of graphs (Section \ref{sec1}, Part I), 
and its different ``variations''  (Part I, \ref{4.20}). 
In the construction, the  resolution graph (string)
of a Hirzebruch-Jung  singularity
plays a distinguished role. For the definition of these strings
and their decorations, see (Part I, \ref{2.3}). 

The algorithm states that the graph $\gxg$ is a cyclic graph--covering   of
$\gc$ with the additional modification which replaces 
each edge by a string (cf. Section \ref{sec1} (Part I) and \ref{4.20}; 
and compare also with
the construction of cyclic coverings in  Section \ref{sec3} and
\ref{2.1}-\ref{2.3} (Part I).  
The point is that for the   covering data of this graph--covering,  one 
can apply proposition (\ref{unique}), which guarantees that up to an
isomorphism, 
there is only one covering graph (with this covering data). 
In particular, the graph--covering is completely determined from the covering 
data, and from the information about the inserted strings. 
Having this in mind, the steps of the algorithm are self--explanatory:

\vspace{1mm}

\noi {\bf Step 1.} provides the covering data of the vertices, and also 
determines some of the (equivariant) decorations of these vertices.

\vspace{1mm}

\noi {\bf Step 2.} provides the covering data of edges, and the description 
of the  inserted strings with their complete decorations.

\vspace{1mm}

\noi {\bf Step 3.} provides the missing decorations of the vertices
constructed in Step 1. After this step we obtain a possible 
embedded resolution graph $\gxg$  of $g|X_k: (X_k,0)\to (\bfc,0)$.

\vspace{1mm}

\noi {\bf Step 4.} is the  procedure of deleting all the arrowheads
 and multiplicities of $\gxg$ (and  successively blowing  down  all 
the rational  $(-1)$--curves) in order to obtain  the (minimal)
 resolution graph of $(X_k,0)$. 

\bekezdes{Definition. -- ``Legs'' and ``stars''.}{ls} 
Before we start the detailed presentation of the algorithm,
we introduce another notion, which helps the presentation  of the
algorithm. 

Fix a  non--arrowhead vertex $v$ of
$\gc$ with weights $(m;n,\nu)$ and $[g]$.
We define the notion of {\em star of $v$}, 
which keeps track of all the edges adjacent to $v$, preserving the 
weights $(n;c,d)$ of the other endpoint of the edges, and unifying the 
different cases represented by loops and 
 edges connecting non--arrowheads or arrowheads. 

The    ``star of  $v$''  consists of the vertex $v$
and different ``legs''. A leg has the following form:

\vspace{1.5cm}

\ \hspace{2cm} \fig{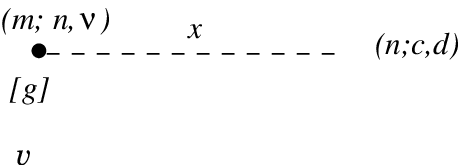}

\noi Here $x\in \{1,2\}$ and the decorations satisfy the same  compatibility 
conditions as the edges (\ref{graphIII}). Sometimes we say that a leg as 
above ``is supported by the vertex $v$''. 

Now, consider again the graph $\gc$ and the fixed non--arrowhead vertex 
$v$. The collection of legs supported by $v$ is defined as follows.

 Any edge, with weight 1 and decoration as in (\ref{summary} Edges, 
Case 1), with one of its endpoints $v_1=v$, provides a leg supported by $v$
with decorations $x=1$ and $(n;c,d)=(m;l,\lambda)$. If the corresponding edge
supports  an arrowhead in $\gc$ then 
automatically  $(n;c,d)=(m;l,\lambda)=(1;0,1)$
(but otherwise, at the level of legs, we do not distinguish between the other 
end of the edge being an arrowhead or not). 

Similarly, any edge, with weight 2 and decoration as in 
(\ref{summary} Edges,  Case 2), with one of its endpoints $v_1=v$, 
provides a leg supported by $v$
with decorations $x=2$ and $(n;c,d)=(m';n,\nu)$. 
Again, if the corresponding edge
supports an arrowhead in $\gc$ then
automatically   $(n;c,d)=(m';n,\nu)=(1;0,1)$.

Any loop supported by $v$ and weighted by $x$, provides  {\em two} 
legs supported by $v$, both decorated by the same $x$ and 
$(n;c,d)=(m;n,\nu)$. 

\vspace{2mm}

The collection of all the legs constructed in this way forms the 
{\em star of $v$}. Obviously, if a non--arrowhead vertex of $\gc$
is connected to $v$ by more than one edge then each edge contributes 
a leg. 
For a general picture of a star,  see the next paragraph.

Another way to define the  star of $v$  is the following. Consider the
topological realization of $\gc$. Then a small neighborhood of 
the point which represents $v$ is the {\em star of $v$}. (Obviously, we have
 to add also the natural decorations.)

\newpage

\subsection*{The algorithm.}

\bekezdes{Step 1. -- The covering data of the vertices.}{s1} \ 

\noi In order to construct the graph $\gxg$ as a  covering graph
of $\gc$  (modified with strings, and decorated as an
embedded resolution graph, see the discussion above) we need the covering data.
The first step provides  the integers $\{n_v\}_{v\in \calv(\gc)}$.

\vspace{1mm}

\noi \alahuz{Case 1.}  \
Consider a non-arrowhead vertex $v$ of \, $\G_\C$ decorated by 
$(m;n,\nu)$ and $[g]$. Consider its star (see the previous paragraph):

\vspace*{2cm}

\fig{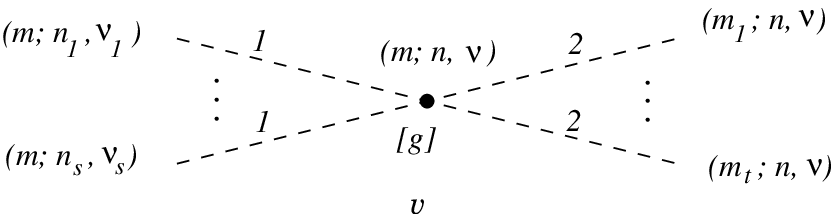}

\noi Let $s$ (respectively $t$) be the number of legs weighted by 
$x=1$ (respectively by   $x=2$). 

\vspace{1mm}

\noi Above the vertex $v$ of $\gc$ put 
$n_v$ non-arrowhead vertices, where  
$$n_v=\gcd(m,k\nu -n,k\nu_1 -n_1,...,k\nu_{s} -n_{s},
m_1,...,m_{t}).$$
Put on each of these non--arrowhead vertices the same decoration:
``the multiplicity'' $(\tilde{m})$:
$$\tilde{m}={m {\nu}\over \gcd (m, {\nu}k-n)}\, ,$$
 and the genus $[\tilde{g}]$
 determined by the formula:

\vspace{2mm}

\begin{tabular}{ll}
$n_v\cdot (2-2\tilde{g}) $ & 
$=(2-2g - s-t) \cdot \gcd (m, k{\nu}-n)$\\
   &$+ \sum\limits_{i=1}^s \gcd (m, k{\nu}-n, k\nu_i- n_i) + 
    \sum\limits_{j=1}^t \gcd (m,k{\nu}- n, m_j).$ \\
\end{tabular}

\noi (In Step 3, the  
``self--intersection'' of each vertex above $v$ will also  be provided.)

\vspace{2mm}

\noi \alahuz{Case 2.} \ Consider an arrowhead vertex $v$ 
\fig{21.eps}  \hspace{1.6cm} of $\G_\C$.
 Above the vertex $v$, in the graph $\gxg$,
put exactly one arrowhead vertex; and let 
 $n_v=1$.  Let the ``multiplicity''  of this 
arrowhead be $1$, i.e. the arrowhead of $\gxg$ is: 
\fig{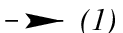}\hspace{1cm}.

\bekezdes{Step 2. -- The covering data of edges and the type of the inserted strings.}{s2} \ 

\vspace{2mm}

\noi This step provides the integers $\{n_e\}_{e\in \cale(\gc)}$. 

\vspace{2mm}

\noi\alahuz{Case 1.}\ Consider an edge $e$ in $\G_\C$, 
 with weight 1  (cf. with \ref{summary} Edges, Case 1):

\vspace*{2cm}

\ \hspace{4cm} \fig{7.eps}

\noi Define:  $$n_e=\gcd(m, k\nu-n, k\lambda-l).$$

\noi 
Notice that Step 1 guarantees that in $\gxg$ there is a column of $n_{v_1}$ 
(resp.  $n_{v_2}$) vertices above the vertex $v_1$ (resp. $v_2$)
of $\G_\C$. Moreover, $n_{v_1}$ and $n_{v_2}$ both divide $n_e$. 

\vspace{1mm}

\noi Then, above the edge $e$   insert (cyclically) in $\gxg$ 
exactly $n_e$ strings of type 
$$Str\left(\, {k\nu-n\over n_e}, {k\lambda-l\over n_e};{m\over n_e}\
\Big|\ \nu, \lambda;0 \, \right).$$

\noi
If the right vertex $v_2$ is replaced by an arrowhead, then the edge $e$ is:

\vspace{1cm}

\ \hspace{4cm} \fig{8.eps}

\noi 
Therefore, in the above  formulae $m=1$, $n_e=1$ 
and the string represents a nonsingular
surface, i.e. it is the empty graph. Therefore, above such an edge $e$
put a single edge (which supports that arrowhead  of $\gxg$ which covers 
the corresponding arrowhead of $\gc$). 

\vspace{2mm}

\noi \alahuz{Case 2.} Consider an  edge $e$  in 
$\G_\C$, with weight $2$  (cf. \ref{summary} Edges, Case 2):

\vspace{2cm}

\ \hspace{4cm} \fig{9.eps}

\noi For such an edge, define:
$$n_e=\gcd(m, m', k\nu-n).$$

\noi 
Similarly as above,
Step 1 guarantees that in $\gxg$ there is a column of $n_{v_1}$ 
(resp.  $n_{v_2}$) vertices above the vertex $v_1$ (resp. $v_2$)
of $\G_\C$, and $n_{v_1}$ and $n_{v_2}$ both divide $n_e$. 

\vspace{1mm}

\noi Then, above the edge $e$   insert (cyclically) in $\gxg$ 
exactly $n_e$ strings of type 
$$Str\left(\,{ m\over n_e},{m'\over n_e};{k\nu-n \over n_e}\
\Big| \ 0,0;\nu\, \right).$$

\noi
If the right vertex $v_2$ is replaced by an arrowhead, then the edge $e$ is:

\vspace{1cm}

\ \hspace{4cm} \fig{10.eps}

\noi 
Therefore, in the above formulae $m'=1$, $n=0$, $\nu=1$ and  $n_e=1$.
In this case, in general,  the string is not empty, it is of type
$Str(m,1;k\,|\, 0,0;1)=Str(m,1;k)$. In particular, above the  edge $e$
(which supports an arrowhead of $\gc$) 
insert a string of type $Str(m,1;k)$
(whose right end  supports that arrowhead  of $\gxg$ which covers 
the corresponding arrowhead of $\gc$). 

\vspace{2mm}

The reader is invited to consult (Part I, \ref{4.20}) and 
 Step 3 of the algorithm presented in (Part I, Section \ref{sec3} and
\ref{2.1}-\ref{2.3}).
where the insertion of the strings is explained in more details. 

\bekezdes{-- The covering data revisited.}{covdata} 
Notice that the  set of integers $\{n_v\}_{v\in\calv(\gc)}$ 
and $\{n_e\}_{e\in\cale(\gc)}$ satisfies the axioms of a covering data.
Moreover,   if $v\in \calv^1(\gc)$ then $m=1$ hence $n_v=1$.
Therefore, by proposition (\ref{unique}), there is only one
cyclic graph--covering  of $\gc$ with this covering data. 
If in this unique graph we replace the edges with the 
corresponding strings (as explained above) we obtain $\gxg$. 

\bekezdes{Step 3. -- The missing decorations.}{s3} The first two steps
provide a graph with some decorations:  all the multiplicities,
all the genera, and some of the self--intersections. 
More precisely, the self--intersections of the vertices
constructed in Step 1 are missing. In this step compute  all these
numbers using  (Part I, \ref{1.3}2). 

\vspace{2mm}

This finishes completely the construction of (a possible) embedded 
resolution graph $\gxg$ of the germ $g|X_k:(X_k,0)\to(\bfc,0)$. 

\bekezdes{Step 4.}{s4} In order  to obtain the resolution graph 
$\gxk$ of $(X_k,0)$,  the arrows and multiplicities
of $\gxg$ have to be deleted. If the resulting  graph is not minimal, 
then blowing down successively
the rational curves with self--intersection $-1$ will provide a minimal one. 

\subsection*{Examples.}

\vspace{2mm}

\bekezdes{Notation.}{notation} In order to help the reader distinguish
the vertices of $\gxg$
 which cover some vertices of $\gc$ (i.e. vertices constructed
in Step 1) from the vertices of the strings inserted in Step 2,
we  indicate the latter by \fig{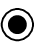} \hspace{1mm}.
Obviously, in the final graph we  make no distinction between them.

In the next examples, we give two graphs, the first is the result of
steps 1-3, the second is the result of step 4. 

\bekezdes{Example.}{e1}
Let $f(x,y,z)=x^2-y^3$ and $g(x,y,z)=z$.
A possible dual graph $\G_\C$ is shown in Example \ref{planecurve}

Set $k=5$. The  algorithm provides the minimal resolution graph
of $(\{x^2+y^3+z^5=0\},0)$.
Applying Steps 1-3, one obtains:

\ \hspace{1cm} \ 

\vspace*{4.5cm}

\fig{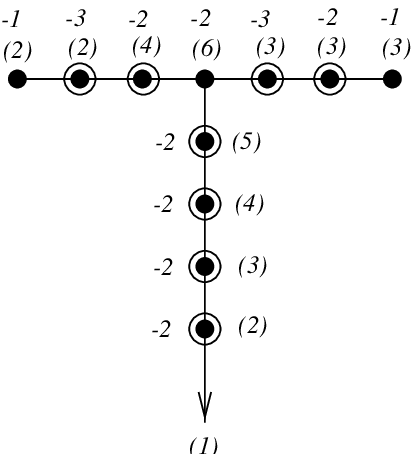} 

\vspace{-3cm} 

\ \hspace{5cm}  Step 4 gives:

\vspace{1cm}

\ \hspace{8cm} \fig{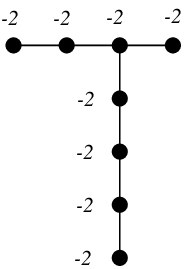}

\vspace{1cm}

\sm\noi This is exactly the expected graph of the surface singularity $E_8$.

\noi
More generally, the above algorithm contains as a particular case
the construction of the resolution graph of surface singularities of type
$(\{f(x,y)+z^N=0\},0)$ (cf. Part I, \ref{8.1}),
via the construction (\ref{planecurve}) of $\gc$. 

\newpage
\bekezdes{Example.}{e2}
Let $f(x,y,z)=x^3+y^2+xyz$ and $g(x,y,z)=z$. 
A possible dual graph $\G_\C$ is shown in Example \ref{andrase}
Let $k=9$, then the algorithm provides the resolution 
of the cusp singularity $(\{x^3+y^2+z^9+xyz\},0)$. 

\vspace{3cm}

\fig{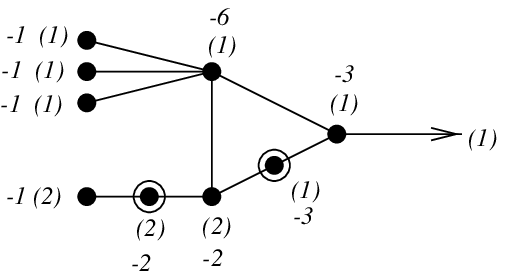}

\vspace{-1cm}

\hspace{8cm} \fig{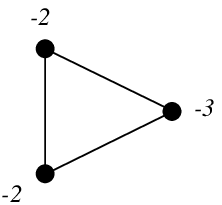}

\vspace{4mm}

\bekezdes{Example.}{e3}
Consider again $f(x,y,z)=x^3+y^2+xyz$ and $g(x,y,z)=z$, as in the previous
example, but take 
$k=12$. Then the algorithm gives the minimal resolution graph
of the cusp singularity $(\{x^3+y^2+z^{12}+xyz\},0)$.

\vspace{3cm}

\fig{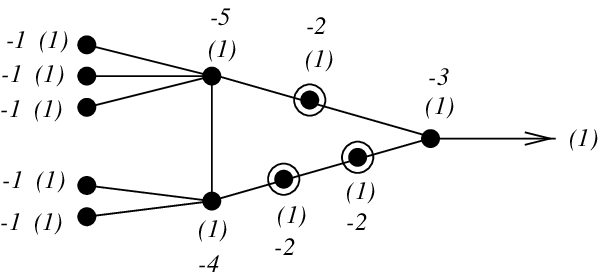}

\vspace{-1cm}

\hspace{8cm} \fig{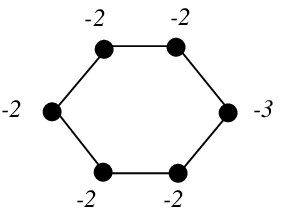}

\vspace{6mm}

\noi 
The interested reader can verify that for arbitrary $k$ (and the same $\Phi$
 and $\gc$ as above),
the graph $\gxk$ is a cyclic graph as above with one $-3$ curve and 
$k-7$ curves with self--intersection $-2$. 

\bekezdes{Example.}{e4}
Let $f(x,y,z)=x^3+y^3+xyz$ and $g(x,y,z)=z$. 
A possible dual graph $\G_\C$ is:

\vspace*{1.5cm}

\begin{tabular}{cc}
\includegraphics{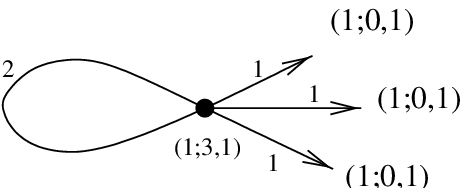}\hspace{6cm}
& (cf.  Example \ref{hurkosn} with $n=3$.)
\end{tabular}

Set $k=6$. Then the minimal resolution graph
of the cusp singularity $(\{x^3+y^3+z^6+xyz\},0)$, obtained in the 
same steps, is:

\vspace*{2.5cm}

\fig{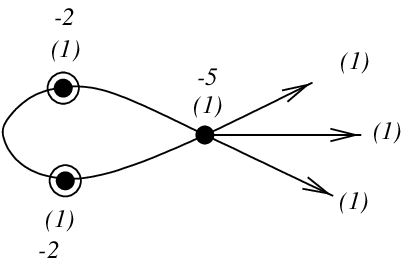}

\vspace{-1cm}

\hspace{8cm} \fig{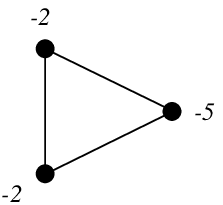}

\vspace{2mm}

\bekezdes{Example.}{e5} 
Consider again $f(x,y,z)=x^3+y^3+xyz$ and $g(x,y,z)=z$ as in
the previous example.  Set $k=8$, then the graph of 
$(\{x^3+y^3+z^8+xyz\},0)$ is:

\ \hspace{2cm} \

\vspace*{2.5cm}

\fig{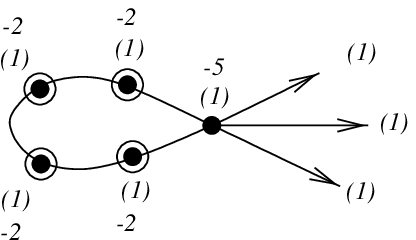}

\vspace{-1cm}

\hspace{8cm}\fig{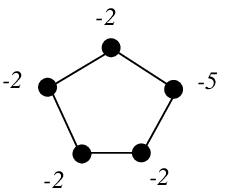}

\vspace{4mm}

\noi 
Again, one  can verify that for arbitrary $k$ (and the same $\Phi$
 and $\gc$ as above),
the graph $\gxk$ is a cyclic graph as above with one $-5$ curve and 
$k-4$ curves with self--intersection $-2$. 

\bekezdes{Example.}{e6}
Let $f(x,y,z)=y^3+(x^2-z^4)^2$ and $g(x,y,z)=z$.
A possible dual graph $\G_\C$ 
is shown in Example (\ref{ujpelda}).
Set $k=15$. The embedded resolution graph of $z$ defined on 
$(\{y^3+(x^2-z^4)^2+z^{15}=0\},0)$ is: 

\vspace{7cm}\hspace{-1cm}
\fig{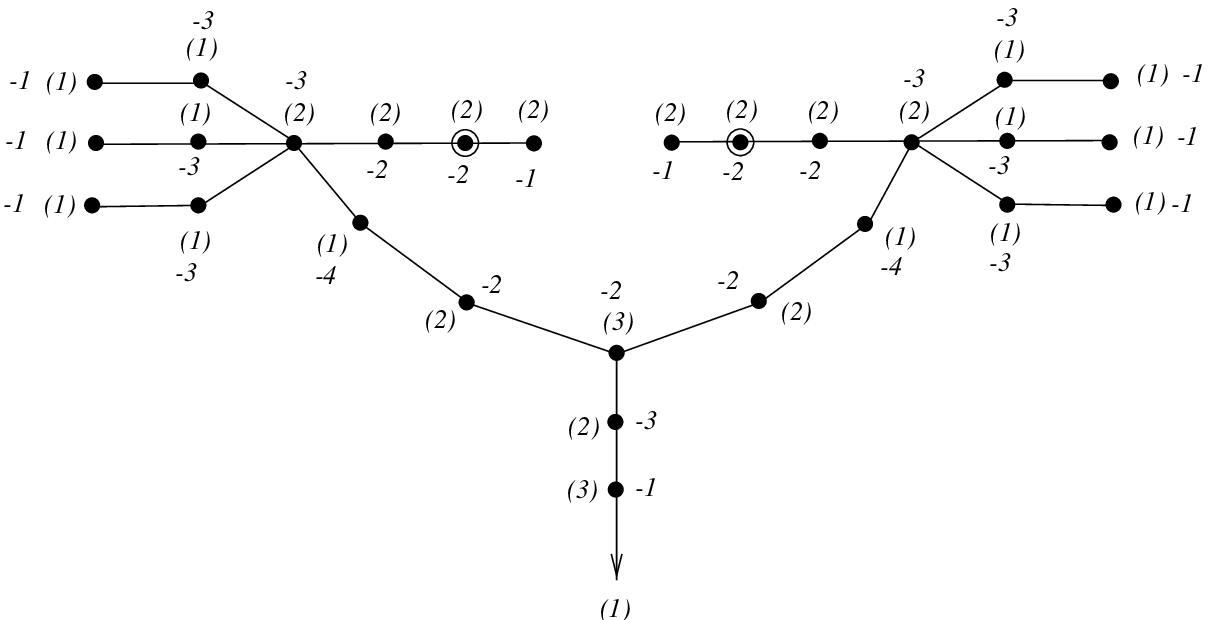}

\noi
Applying Step 4, one gets the resolution graph of
$(\{y^3+(x^2-z^4)^2+z^{15}=0\},0)$: 

\vspace*{4cm}

\hspace{2cm} \fig{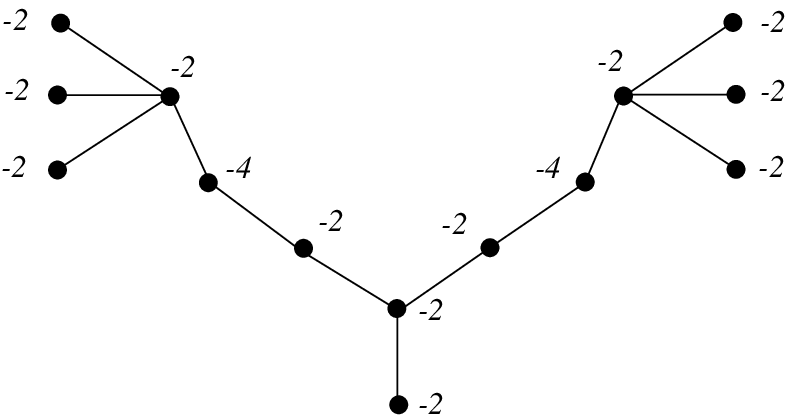}

\bs\noi
Note that this result can be checked by using
the algorithm described in (Part I, \ref{8.1}), with
$f(x,y)=(x^2-y^4)^2+y^{15}$ and $N=3$.

\newpage

\subsection*{Proof of the algorithm and final remarks.}

\vspace{2mm}

\bekezdes{-- Outline of the proof.}{hint} 
In this subsection we present the main steps of  the proof of the algorithm, 
we emphasize 
the key points, and leave all  the details to the reader.

\vspace{2mm}

The  fundamental fact is that the lift of the singularity $(X_k,0)$ in
$V^{emb}$, via the embedded resolution $r: (V^{emb},\de)\to (\bfc^3,D)$,
has  a distinguished position. Namely, let the strict transform of $X_k$
be denoted by $\xst$. Then for  sufficiently large $k$
$$ \xst\cap \de=\C \ \ \ \ \mbox{and} \ \ \ Sing(\xst)\subset \C.$$
Moreover, along an irreducible component of the intersection 
$\xst\cap\dek$, the surface 
$\xst$ is either smooth or forms an equisingular family of plane
curve singularities. Therefore, 
if $\xst$ is normalized, the  singularities 
along the intersection $\xst\cap\dek$ are smoothed.
The resulting  normal surface $\widetilde{\xst}$ has 
singularities only above  points of $\xst\cap \deh$.
Moreover, the singularities of  the strict transform 
$\xst$, above the points  $\xst\cap\deh$, 
have local equations of type  $\{ w^\alpha = u^\beta v^\gamma\}$, with
$\beta, \gamma\geq 1$ and $\alpha >1$ (in order to find these integers,
use the local equations (\ref{summary})). 
Thus, after normalization, 
these singularities are of Hirzebruch-Jung type  and
their resolution is well-known (cf. Part I, \ref{2.3}). 
A smooth surface $X_k^{res}$, which
provides a resolution of the original isolated surface singularity
$(X_k,0)$, can be obtained by resolving
these Hirzebruch--Jung singularities.

\vspace{2mm}

The construction is summarized in the next diagram:
     
$$ X_k^{res}\ \stackrel{R}{\longrightarrow}\ 
\widetilde{\xst}\ \stackrel{n}{\longrightarrow}\
\xst \ \stackrel{r|}{\longrightarrow}\ (X_k,0),$$

\noi where $r|$ is the restriction of the fixed embedded resolution
$r$; $n$ is the normalization
of $\xst$ and $R$ is the map resolving the 
Hirzebruch-Jung singularities of the surface $\nxst$.

The graph constructed using the algorithm of  the previous section is exactly 
the dual resolution graph of the resolution $(X_k^{res}, \ r|\circ n\circ R)$.
Note that by the above discussion the exceptional curves of the
resolution $X^{res}_k$ are of two type.
Some exceptional curves are irreducible components of the lift
$n^{-1}(\C)$ and these are lifted further isomorphically by $R$.
Other exceptional curves 
appear as a result of resolving Hirzebruch-Jung strings. 

\vspace{2mm}

The crucial fact in the proof is that
 all the information necessary to recover the dual 
embedded resolution graph can be obtained by local multiplicity (vanishing 
order) computations. This claim needs some explanation. 

\vspace{1mm}

The first technical difficulty (the ``easy one'') arises because of the fact
 that the dual resolution graph includes 
the self-intersection numbers of the exceptional  curves as weights.
These numbers are global invariants (characteristic
classes) and in general, their computation is difficult.
But if we consider the germ of an analytic function defined on $X_k$,
then in any embedded resolution graph of it,  
all the self--intersections can be 
computed  from the multiplicities of the germ using  (Part I, \ref{1.3}2). 
This  is exactly  the motivation to consider the {\em embedded resolution }
of the restriction of $g$ to $(X_k,0)$ instead of the resolution of $(X_k,0)$.
It is not difficult 
to verify that the procedure outlined above in fact gives the 
embedded resolution of  $g|X_k:(X_k,0)\to (\bfc,0)$, and 
all the vanishing orders can be determined by local computations. 

\vspace{2mm}

The second (globalization) problem is more serious. First we would like to 
remind the reader about the result regarding the resolution graph of 
cyclic covering. Fix a normal surface singularity $(X,0)$ and a germ
$f$ defined on $(X,0)$, and consider the $N$--cyclic covering $X_{f,N}$ of
$(X,0)$ branched along $\{f=0\}$. Then in general, it is impossible to 
recover the resolution graph of $X_{f,N}$ from the embedded resolution graph
of $f$ and the integer $N$. Even if the embedded resolution graph of $f$ and 
the integer $N$ carry all the needed {\em local} information about the covering
and multiplicities, the {\em global}  information about monodromies,
which determines  the type of the graph--covering (i.e. the covering
data) and also the shape of the covering graph, is missing. 
In this case of cyclic coverings, 
this global information can be codified in the universal covering graph 
of the embedded resolution graph of $f$. 
On the other hand, if the embedded resolution graph of $f$ has some additional
 properties (e.g. if it is a tree and $g_w=0$ for all vertices $w\in\calw$),
then the necessary  global information can be recovered from the local one. 
(For more details,  see Part I of the paper.)

\vspace{1mm}

In the case of the generalized Iomdin series 
 we can ask the same  question:  is the global information 
(i.e. the number of irreducible components of $n^{-1}(C)$ in 
$\widetilde{\xst}$ situated above an irreducible  component $C$  of $\C$ in
$\xst$, or the global configuration of the intersection of these components)
codified in the local data? The answer is  {\em yes}!  This follows 
from  the results of Section \ref{sec5}, 
where it is proven 
that the covering data can be determined from $\gc$.
Moreover, proposition (\ref{unique})
(cf. also Part I, \ref{theo2})  guarantees that with this covering data
there is only one graph--covering. 
The fact which 
makes these claims true is that that part of $\C$ (namely $\C^2$),  
above which the covering is not
 trivial, contains only rational curves, and the dual graph of these curves 
has  no cycles.

\vspace{2mm}

Once it is 
clarified that no global information is needed, everything is reduced 
to local  computations. They  are left to the reader. \hfill $\diamondsuit$

\bekezdes{Final remarks.}{fr} \ 

{\bf 1.}\ If $F_h$ denotes the Milnor fiber of a hypersurface singularity
$h:(\bfc^3,0)\to (\bfc,0)$, and $\chi$  the topological Euler--characteristic,
then one has the following formula:
$$\chi(F_{f+g^k})-\chi(F_f)=k\cdot \,\sum_{j=1}^s\, d_j\mu_j,$$
where $\mu_j$ denoted the Milnor number of the transversal singularity
$T\Sigma_j$ associated with $\Sigma_j$ ($1\leq j\leq s$). This
formula for the ``classical 
case'' (when $g$ is a generic linear form) was proved by Iomdin 
\cite{Io}, for the general case, see \cite{NZ,NP}.

Notice that both integers $d_j$ and $\nu_j$ can be determined from the 
weighted graph  $\gc$. Indeed, $d_j$ is given in (\ref{dj}). On the other 
hand, theorems (\ref{cts1}) and (\ref{cts2}) provide the embedded
resolution graph of $T\Sigma_j$ which via a result of A'Campo's 
\cite{A'2} provides $\mu_j$. 

In particular, $\gc$ determines completely the correction term 
$i(f+g^k)-i(f)$ in the case of $i(h)=\chi(F_h)$ (cf. \ref{invs}). 

\vspace{2mm}

{\bf 2.}\ Notice that for any vertex $w\in \calw(\gce)$ the covering data 
is $n_w=1$, and also for any edge $e\in \cale(\gce)$ one has $n_e=1$. 
Moreover, 
the type of the inserted string corresponding to $e$ is trivial.
Therefore, above $\gce$ the covering is trivial, i.e. for any $k$, the graph 
(above) $\gce$ stays as a stable part of the resulting graph $\gxg$ or $\gxk$. 

On the other hand, 
the graph above $\gck$ is increasing as $k$ becomes larger. 
More precisely, the number of vertices above a vertex $w\in \calw(\gck)$
is a periodical function in $k$.
Similarly,  the 1--edges of $\gck$ are replaced by strings which behave
periodically  when $k$ runs over the integers $k\geq k_0$. But the 2-edges
of $\gck$ are replaced by strings whose ``size''  increases infinitely as
$k\to\infty$. Basically, the union of these 
 increasing spots, associated with the 
subgraphs $\gj$ ($1\leq j\leq s$), correspond to the ``correction term''
$i(f+g^k)-i(f)$, where $i$ denotes the resolution graph (cf. also 
with \ref{invs}). 

\vspace{2mm}

{\bf 3.}\ Not all the decorations of the weighted graph $\gc$ were used 
in the above algorithm. Nevertheless, we believe that all of them have 
important geometrical significance (we will explain this  in the forthcoming
 paper \cite{NSz}, where we plan to discuss the relationship of $\gc$ with 
other invariants). 

\vspace{2mm}

{\bf 4.}\ Any fault--finder reader can object the following. We constructed
our resolution graph $\gxg$ from the graph $\gc$ whose construction is based
on some properties of the embedded resolution $r$, which is a more complicated
object than any resolution of $(X_k,0)$. This is perfectly true!

But, if we consider our results as qualitative results, then already the 
existence of a graph $\gc$ which coordinates all the resolutions graphs
$\{\gxg\}_k$ in a purely combinatorial way (as it is explained in the 
algorithm), and has all the properties listed in Section \ref{sec5},
 (and 
conjecturally  provides all the ``correction terms'' $i(f+g^k)-i(f)$
for different invariants $i$ such as the zeta function, signature, spectral 
pairs, etc.), is really remarkable. 

\vspace{2mm}

{\bf 5.}\ In this paper  the case of the Iomdin series was  discussed, 
but all the results can be generalized to the case of an
 arbitrary series  of composed singularities.

{

}

\end{document}